\documentclass[11pt,a4paper]{article}
\usepackage[T1]{fontenc}
\usepackage[utf8]{inputenc}
\usepackage{lmodern}
\usepackage[margin=1in]{geometry}
\usepackage{amsmath,amssymb,amsthm,mathtools}
\usepackage{booktabs}
\usepackage{graphicx}
\usepackage[section]{placeins}
\usepackage{microtype}
\usepackage[round,authoryear]{natbib}
\usepackage[colorlinks,linkcolor=black,citecolor=black,urlcolor=black]{hyperref}

\newtheorem{theorem}{Theorem}[section]
\newtheorem{proposition}[theorem]{Proposition}
\newtheorem{lemma}[theorem]{Lemma}
\newtheorem{corollary}[theorem]{Corollary}
\theoremstyle{definition}
\newtheorem{definition}[theorem]{Definition}
\newtheorem{example}[theorem]{Example}
\theoremstyle{remark}
\newtheorem{remark}[theorem]{Remark}

\newcommand{\E}{\mathbb{E}}
\renewcommand{\P}{\mathbb{P}}
\newcommand{\Var}{\operatorname{Var}}

\newcommand{\cx}{\mathrm{cx}}
\newcommand{\dd}{\mathrm{d}}
\newcommand{\Cr}{C^{\uparrow}}
\newcommand{\Rrho}{R_\rho}
\newcommand{\Fp}{F_\phi}
\newcommand{\norm}[1]{\lVert#1\rVert}
\newcommand{\sprod}[2]{\langle#1,#2\rangle}

\title{Dependence measures as spectral risk measures}
\author{Karl Friedrich Siburg\thanks{Fakult\"at f\"ur Mathematik, Technische Universit\"at Dortmund, Germany.}}
\date{}

\begin{document}
\maketitle

\begin{abstract}
\noindent
For a bivariate copula $C$ with $(U,V)\sim C$, the conditional distribution functions $L_t=\partial_1C(U,t)$ are random variables with means $t$, so that measures of directed dependence become risk functionals. The absolute Lorenz curves of the variables $L_t$ form the rearranged copula $C^\uparrow$. In this picture Chatterjee's rank correlation $\xi$ is a squared Sobolev seminorm and the rearranged Spearman correlation an $L^1$ norm of $C^\uparrow-\Pi$, and both arise from one spectral risk measure, applied either to the conditional laws of $V$ or to the laws of $L_t$. The resulting classes $\omega_\phi$ and $\kappa_\phi$ of spectral dependence measures have mirror-image axiomatic properties. For every spectrum we determine the exact upper bound of $\kappa_\phi$ in terms of $\xi$, and we prove sharp inequalities between $\xi$, $\kappa_\phi$ and the values of $C^\uparrow$, which show how little global measures say about the tails. For plug-in estimators based on the empirical checkerboard copula we obtain consistency, rates and a resolution limit, and a simulation study shows that tail-localised measures separate copulas with identical $\xi$. An application to large credit portfolios illustrates the results.


\end{abstract}

\section{Introduction}

A measure of (directed) dependence assigns to a pair of continuous random variables $(X,Y)$ a number in $[0,1]$ that vanishes exactly when $X$ and $Y$ are independent and equals one exactly when $Y$ is a measurable function of $X$. The measure $\xi$ was introduced for continuous random variables by \citet{DSS2013}; up to an affine transformation it is the squared $L^2$ norm of the partial derivative $\partial_1C$ of the copula $C$, that is, the part of the squared Sobolev norm of \citet{SiburgStoimenov2008} that belongs to the conditioning variable. It became widely known through the simple rank estimator of \citet{Chatterjee2021}, see also \citet{AzadkiaChatterjee2021}, and is now often called Chatterjee's rank correlation. Further examples are Trutschnig's $\zeta_1$ \citep{Trutschnig2011}, the rearranged dependence measures of \citet{SDS2024}, and the sensitivity measures of \citet{ALFT2026}. All of them are functionals of the conditional distribution functions $u\mapsto\partial_1C(u,t)$, and all of them are monotone in the Schur order for conditional distributions \citep{AnsariRueschendorf2021,AnsariRockel2024,AnsariFuchs2025}.

This paper develops a simple idea: such measures are \emph{risk functionals}. If $(U,V)\sim C$, then $L_t:=\partial_1C(U,t)=\P(V\le t\mid U)$ is a random variable with values in $[0,1]$ and mean $t$, and the Schur order compares the families $(L_t)_t$ in convex order. Dependence is therefore riskiness of the variables $L_t$ at fixed mean. We use two classical diagrams from the theory of stochastic dominance, the absolute Lorenz curve and the stop-loss transform, which are Fenchel conjugates of each other \citep{OgryczakRuszczynski2002}, together with spectral risk measures \citep{Acerbi2002} and the dual theory of choice of \citet{Yaari1987}.

\paragraph{Main results.}
\emph{(i) Two diagrams.} The Lorenz curves of the variables $L_t$ form the rearranged copula $\Cr$ of \citet{SDS2024}, and Fenchel duality gives
\[
\Cr(x,t)=\min_{q\in[0,1]}\Bigl\{xq+\int_0^1\bigl(\partial_1C(u,t)-q\bigr)_+\,\dd u\Bigr\},
\]
which determines $\Cr$ without any rearrangement (Lemma~\ref{lem:basic}).

\emph{(ii) A dictionary.} With $h=\Cr-\Pi$ we have $\xi(C)=6\norm{\partial_1h}_{L^2}^2$, $\Rrho(C)=12\norm{h}_{L^1}$ and $\zeta_1(C)=6\int_0^1\norm{h(\cdot,t)}_\infty\,\dd t$ (Theorem~\ref{thm:dictionary}); in particular $\Rrho=\Lambda_1$ with the $L^1$ sensitivity measure of \citet{ALFT2026}. We obtain $\tfrac23\zeta_1^2\le\xi\le\zeta_1\le\Rrho\le\min\{2\zeta_1-\zeta_1^2,\sqrt{2\xi}\}$ with optimal constants and the exact region of the pairs $(\zeta_1,\Rrho)$ (Corollary~\ref{cor:ineq}).

\emph{(iii) Spectral dependence measures.} The spectrum $\phi_1(s)=2(1-s)$, applied to the conditional laws of $V$, gives $\xi$, and applied to the laws of $L_t$ it gives $\Rrho$ (Proposition~\ref{prop:two}). A general spectrum $\phi$ gives two classes $\omega_\phi$ and $\kappa_\phi$, related by Yaari's exchange of outcomes and probabilities. The measures $\kappa_\phi$ are the scalar products of $\Cr-\Pi$ with fixed directions in the Hilbert space structure in which $\xi$ is the squared norm (Proposition~\ref{prop:pairing}). For every spectrum we determine the exact upper bound of $\kappa_\phi$ in terms of $\xi$ together with the extremal copulas; for $\phi_1$ this is the upper boundary of the $(\xi,\rho)$ region of \citet{AnsariRockel2025} (Theorem~\ref{thm:exact}). The same method gives the exact joint range of two measures $\kappa_\phi$ and $\kappa_\psi$ (Proposition~\ref{prop:pairs}). The axiomatic properties are mirror images: $\omega_\phi=1$ characterises complete dependence and $\kappa_\phi=0$ independence, whereas each of the two remaining characterisations holds if and only if $\phi$ is strictly decreasing (Theorem~\ref{thm:mirror}). Near independence $\kappa_\phi$ is of first and $\omega_\phi$ of at most second order (Proposition~\ref{prop:order}).

\emph{(iv) Sharp inequalities.} For every spectrum and every point $(x,t)$ we determine the smallest values of $\xi(C)$ and $\kappa_\phi(C)$ compatible with a given value of $\Cr(x,t)$; the extremal copulas are explicit and the same for all spectra (Theorem~\ref{thm:sharp}). Consequences are sharp bounds for $\norm{C-\Pi}_\infty$ in terms of $\xi$ and $\Rrho$ together with a sharp reverse bound (Corollary~\ref{cor:uniform}), and the fact that no measure with a bounded spectrum controls $\Cr$ near the boundary of the unit square (Corollary~\ref{cor:tail}), which tail-localised spectra repair (Corollary~\ref{cor:stress}). Analogous results hold for the stop-loss diagram (Proposition~\ref{prop:profile} and Theorem~\ref{thm:stoploss}).

\emph{(v) Estimation.} Plug-in estimators based on the empirical checkerboard copula are strongly consistent and, under the Lipschitz condition of \citet{AnsariWiesel2026}, converge at the rate $N^{-1/3}$, which cannot be improved for $\omega_\phi$ (Theorem~\ref{thm:estimation}); tail levels below $N^{-1/3}$ cannot be resolved (Proposition~\ref{prop:resolution}). In a simulation study, tail-localised measures separate a Gaussian and a Student copula with equal $\xi$ from a few thousand observations, whereas all global measures fail to do so (Section~\ref{sec:simulation}).

\emph{(vi) An application.} In factor models of portfolio credit risk, $L_t$ is the loss of a large homogeneous portfolio with default probability $t$. Then $\xi$ is the non-diversifiable share of the default variance, the two diagrams give expected shortfall and tranche losses, and the sharp inequalities determine what a single dependence parameter can say about tail risk (Section~\ref{sec:credit}).

\paragraph{Related work.}
The inequality for $\xi$ in (iv) extends the bound $|\beta|^3\le2\xi$ of \citet{OrendayRockel2026} for Blomqvist's $\beta$, which is the case $x=t=\tfrac12$. The inequality $\xi\le\rho$ for stochastically increasing copulas contained in (ii) is due to \citet{AnsariRockel2025}, who also determine the exact $(\xi,\rho)$ region; applied to $\Cr$, their result gives the exact upper bound for $\Rrho$ in terms of $\xi$, of which our bound $\sqrt{2\xi}$ is the asymptote at the origin (Remark~\ref{rem:region}). Theorem~\ref{thm:exact} extends this upper bound to all spectral measures $\kappa_\phi$ and gives a short proof for it. Integral functionals of the form $\int\!\!\int\Psi(\partial_1C(u,t)-t)$ and $\int\!\!\int\!\!\int\Psi(\partial_1C(u,t)-\partial_1C(u',t))$ with convex $\Psi$ are studied in \citet{AnsariFuchs2025} and \citet{ALFT2026}, respectively; our measures $\omega_\phi$ are of a different, level-dependent form and coincide with these only for affine spectra, all of which give Chatterjee's $\xi$ (Remark~\ref{rem:lit}). A general theory of measures of the form $\E[\zeta(\P^Y,\P^{Y\mid X})]$ with convex $\zeta$, including necessary and sufficient conditions for the characterisation of independence and of functional dependence, is due to \citet{BorgonovoEtAl2025}; the measures $\omega_\phi$ are of this form, and the corresponding half of Theorem~\ref{thm:mirror} is an instance of their results (Remark~\ref{rem:lit}). Standardised linear functionals of the rearranged copula appear as examples of rearranged dependence measures in \citet{AnsariWiesel2026}; the measures $\kappa_\phi$ are of this type, and the convergence rates for checkerboard estimators obtained there apply to both of our classes (Section~\ref{sec:estimation}). A different link between rearranged measures and other classes of measures is established in \citet{Ansari2026}, and the representation of $\xi$ through conditionally independent copies is also the starting point of \citet{Limbach2026}. Gini-type risk functionals related to our spectra appear in \citet{Yitzhaki1982}, \citet{OgryczakRuszczynski2002} and \citet{FurmanWangZitikis2017}. The interpretation of $L_t$ as the loss of a large credit portfolio in Section~\ref{sec:credit} connects our results with \citet{AnsariLuetkebohmert2025}, who compare Bernoulli mixture models in convex order.

\section{Preliminaries}\label{sec:prelim}

\subsection{Conditional distribution functions as random variables}

Let $C$ be a bivariate copula, $(U,V)\sim C$, and let $K_C$ be the Markov kernel of $C$, so that $K_C(u,[0,t])=\partial_1C(u,t)$ for almost every $u$; we fix a version for which $t\mapsto K_C(u,[0,t])$ is a distribution function for every $u$. We write $\Pi(u,v)=uv$ and $M(u,v)=\min\{u,v\}$, and we call $C$ \emph{completely dependent} if $V=f(U)$ almost surely for a measurable $f$. For $t\in[0,1]$ put
\[
L_t:=K_C(U,[0,t])=\P(V\le t\mid U).
\]
Then $0\le L_t\le1$, $\E L_t=t$, and $t\mapsto L_t$ is nondecreasing, so that $\E|L_{t'}-L_t|=|t'-t|$. We have $L_t\equiv t$ for all $t$ if and only if $C=\Pi$, and $L_t\in\{0,1\}$ almost surely for all $t$ if and only if $C$ is completely dependent. When necessary we write $L^C_t$. A copula is \emph{stochastically increasing} if $u\mapsto\partial_1C(u,t)$ is nonincreasing for every $t$.

For a random variable $L$ with values in $[0,1]$ let $q_L$ be its quantile function and $L^*(s):=q_L(1-s)$ its decreasing rearrangement. We recall that
\begin{equation}\label{eq:supform}
\int_0^xL^*(s)\,\dd s=\sup\bigl\{\E[L\eta]:\ \eta\text{ random variable},\ 0\le\eta\le1,\ \E\eta=x\bigr\}.
\end{equation}

\begin{definition}\label{def:diagrams}
For $t\in[0,1]$ and $x,q\in[0,1]$ let
\begin{align*}
S_t(x)&:=\int_0^xL_t^*(s)\,\dd s, & h_t(x)&:=S_t(x)-xt,\\
\Theta_t(q)&:=\E(L_t-q)_+, & \theta_t(q)&:=\Theta_t(q)-(t-q)_+ .
\end{align*}
\end{definition}

$S_t$ is the absolute Lorenz curve of $L_t$, taken from the top, and $\Theta_t$ is its stop-loss transform. The functions $h_t$ and $\theta_t$ are the excesses of the two diagrams over their values under independence.

\begin{lemma}\label{lem:basic}
Let $C$ be a copula and $x,q,t\in[0,1]$.
\begin{enumerate}
\item[\textup{(a)}] $S_t(x)=\Cr(x,t)$, where $\Cr$ is the rearranged copula of \citet{SDS2024}. Thus $h_t(x)=(\Cr-\Pi)(x,t)$.
\item[\textup{(b)}] $h_t$ is concave with $h_t(0)=h_t(1)=0$, and $0\le h_t(x)\le M(x,t)-xt$.
\item[\textup{(c)}] For $t<t'$ we have $0\le S_{t'}(x)-S_t(x)\le t'-t$. Hence $t\mapsto h_t(x)$ is Lipschitz continuous with slopes in $[-x,1-x]$.
\item[\textup{(d)}] $0\le\theta_t(q)\le M(q,t)-qt$.
\item[\textup{(e)}] \textup{(Duality.)} $S_t(x)=\min_{q\in[0,1]}\{xq+\Theta_t(q)\}$ and $\Theta_t(q)=\max_{x\in[0,1]}\{S_t(x)-xq\}$. In particular $\max_xh_t(x)=\theta_t(t)=\E(L_t-t)_+=\tfrac12\E|L_t-t|$.
\item[\textup{(f)}] $|C(x,t)-xt|\le\max\{h_t(x),h_t(1-x)\}$.
\item[\textup{(g)}] The map $q\mapsto\theta_t(q)$ has slopes in $[0,1]$ on $[0,t]$ and in $[-1,0]$ on $[t,1]$; the map $t\mapsto\theta_t(q)$ has slopes in $[0,1]$ on $[0,q]$ and in $[-1,0]$ on $[q,1]$.
\end{enumerate}
\end{lemma}

\begin{proof}
(a) This is the definition of $\Cr$ in \citet[Definition~2.2]{SDS2024}, who rearrange $u\mapsto\partial_1C(u,t)$ decreasingly and integrate. (b) $S_t$ is the integral of a nonincreasing function and $S_t(1)=t$. The choice $\eta\equiv x$ in \eqref{eq:supform} gives $S_t(x)\ge xt$, and $0\le L_t\le1$ gives $S_t(x)\le\min\{x,S_t(1)\}$. (c) The first inequality holds since $L_{t'}\ge L_t$. If $\eta$ is optimal in \eqref{eq:supform} for $S_{t'}(x)$ (an optimal $\eta$ exists, for instance $\mathbf 1\{L_{t'}>s\}+\gamma\mathbf 1\{L_{t'}=s\}$ with $s=L_{t'}^*(x)$ and a suitable $\gamma\in[0,1]$), then $S_{t'}(x)=\E[L_t\eta]+\E[(L_{t'}-L_t)\eta]\le S_t(x)+t'-t$. (d) Nonnegativity is Jensen's inequality, and $(\ell-q)_+\le(1-q)\ell$ for $\ell\in[0,1]$ gives $\Theta_t(q)\le(1-q)t$, which is the upper bound in both cases $q\ge t$ and $q<t$. (e) The first formula is the representation of expected shortfall by \citet{RockafellarUryasev2002}; the second follows by Fenchel--Moreau since $S_t$ is concave and continuous, see also \citet[Theorem~3.1]{OgryczakRuszczynski2002}. For $q=t$ note that $\E(L_t-t)_+=\E(t-L_t)_+$. (f) By \eqref{eq:supform}, $C(x,t)=\E[L_t\mathbf 1\{U\le x\}]\le S_t(x)$ and $C(x,t)=t-\E[L_t\mathbf 1\{U>x\}]\ge t-S_t(1-x)$. (g) $\Theta_t$ is convex in $q$ with slope $-\P(L_t>q)\in[-1,0]$, and $(t-q)_+$ has slope $-1$ on $[0,t)$ and $0$ on $(t,1]$. For fixed $q$ and $t<t'$, monotonicity of $t\mapsto L_t$ gives $0\le\Theta_{t'}(q)-\Theta_t(q)\le\E(L_{t'}-L_t)=t'-t$, and $t\mapsto(t-q)_+$ has slope $0$ for $t<q$ and $1$ for $t>q$.
\end{proof}

Part (f) contains the inequalities $C^\downarrow\le C\le\Cr$ of \citet[Lemma~A.1]{SDS2024}, since $t-S_t(1-x)$ is the decreasingly rearranged copula $C^\downarrow(x,t)$.

\subsection{Dependence measures and the Schur order}

We consider
\begin{align*}
\xi(C)&=6\int_0^1\!\!\int_0^1\bigl(\partial_1C(u,t)\bigr)^2\,\dd u\,\dd t-2, &
\Rrho(C)&=12\int_0^1\!\!\int_0^1\Cr(x,t)\,\dd x\,\dd t-3,\\
\zeta_1(C)&=3\int_0^1\!\!\int_0^1\bigl|\partial_1C(u,t)-t\bigr|\,\dd u\,\dd t, &
\Lambda_1(C)&=3\int_0^1\!\!\int_0^1\!\!\int_0^1\bigl|\partial_1C(u,t)-\partial_1C(u',t)\bigr|\,\dd u\,\dd u'\,\dd t,
\end{align*}
that is, the measure of \citet{DSS2013} (Chatterjee's rank correlation), the rearranged Spearman correlation $\rho(\Cr)$ of \citet{SDS2024}, the measure of \citet{Trutschnig2011}, and the $L^1$ member of the family of \citet{ALFT2026}. We also use Trutschnig's metric $D_1(C,D)=\int_0^1\!\int_0^1|\partial_1C-\partial_1D|\,\dd u\,\dd t$, so that $\zeta_1(C)=3D_1(C,\Pi)$.

On the space of functions $A$ on $[0,1]^2$ for which $A(\cdot,t)$ is absolutely continuous for almost every $t$ and $\partial_1A\in L^2([0,1]^2)$, which contains the copulas, we consider the bilinear form and seminorm
\[
\sprod AB_1:=\int_0^1\!\!\int_0^1\partial_1A\,\partial_1B\,\dd u\,\dd t,\qquad\norm{A}_1:=\sprod AA_1^{1/2},
\]
that is, the part of the Sobolev scalar product of \citet{SiburgStoimenov2008} that belongs to the conditioning variable. Since $\sprod C\Pi_1=\sprod\Pi\Pi_1=\tfrac13$ for every copula $C$, we have $\xi(C)=6\norm C_1^2-2=6\norm{C-\Pi}_1^2$.

The Schur order for conditional distributions reads, in our notation,
\begin{equation}\label{eq:schur}
C\le_SD\quad:\Longleftrightarrow\quad L^C_t\le_{\cx}L^D_t\ \text{ for all }t
\quad\Longleftrightarrow\quad\Cr\le D^\uparrow\ \text{ pointwise},
\end{equation}
where $\le_{\cx}$ denotes the convex order and the equivalence is the Hardy--Littlewood--P\'olya characterisation for variables with equal means \citep[Section~3.A]{ShakedShanthikumar2007}; see \citet{AnsariRueschendorf2021} and \citet{SDS2024}. The independence copula is the smallest element and the completely dependent copulas are the largest elements. By \eqref{eq:schur}, every family of functionals of the laws of $L_t$, $t\in[0,1]$, that are consistent with the convex order defines, after integration over $t$, a functional of $C$ that is monotone in $\le_S$.

\section{Norms and areas}\label{sec:dictionary}

Let $V'$ be a copy of $V$ that is conditionally independent of $V$ given $U$. Then $\P(V\le t,V'\le t)=\E L_t^2$, which is the diagonal of the copula of $(V,V')$, the Markov product of the transpose of $C$ with $C$. For $t\in(0,1)$ we call
\[
\varrho_C(t):=\frac{\Var(L_t)}{t(1-t)}=\frac{\P(V\le t,V'\le t)-t^2}{t(1-t)}\in[0,1]
\]
the \emph{indicator correlation} of $C$ at level $t$; it is the correlation of the indicators of $\{V\le t\}$ and $\{V'\le t\}$.

\begin{theorem}[Dictionary]\label{thm:dictionary}
Let $C$ be a copula, $h=\Cr-\Pi$, and let $L_t'$ be an independent copy of $L_t$.
\begin{enumerate}
\item[\textup{(i)}] $\displaystyle\xi(C)=6\int_0^1\Var(L_t)\,\dd t=\int_0^1\varrho_C(t)\,6t(1-t)\,\dd t$, and\\[0.5ex]
$\displaystyle\xi(C)=6\norm{h}_1^2=6\int_0^1\!\!\int_0^1\bigl(\partial_xh_t(x)\bigr)^2\,\dd x\,\dd t=12\int_0^1\!\!\int_0^1\theta_t(q)\,\dd q\,\dd t$.
\item[\textup{(ii)}] $\displaystyle\Rrho(C)=12\int_0^1\!\!\int_0^1h_t(x)\,\dd x\,\dd t=3\int_0^1\E|L_t-L_t'|\,\dd t=\Lambda_1(C)$.
\item[\textup{(iii)}] $\displaystyle\zeta_1(C)=6\int_0^1\max_{x}h_t(x)\,\dd t=6\int_0^1\theta_t(t)\,\dd t$.
\end{enumerate}
\end{theorem}

\begin{proof}
(i) Since $\int\!\!\int(\partial_1C)^2=\int_0^1\E L_t^2\,\dd t=\int_0^1\Var(L_t)\,\dd t+\tfrac13$, the first line holds by the definition of $\varrho_C$. The norm identity follows from $\partial_xh_t(x)=L_t^*(x)-t$. For the last expression, Fubini's theorem gives $\int_0^1\Theta_t(q)\,\dd q=\tfrac12\E L_t^2$ and $\int_0^1(t-q)_+\,\dd q=\tfrac12t^2$, hence
\begin{equation}\label{eq:areaPrimal}
\int_0^1\theta_t(q)\,\dd q=\tfrac12\Var(L_t).
\end{equation}
(ii) By Lemma~\ref{lem:basic}(a), $\Rrho(C)=12\int\!\!\int(\Cr-\Pi)$. For fixed $t$, with $L=L_t$,
\begin{equation}\label{eq:areaDual}
\int_0^1h_t(x)\,\dd x=\int_0^1(1-s)\bigl(L^*(s)-t\bigr)\,\dd s=\tfrac12\int_0^1(2r-1)\,q_L(r)\,\dd r=\tfrac14\E|L-L'|,
\end{equation}
the last step being the representation of Gini's mean difference as an $L$-statistic. The right-hand side of (ii) is the definition of $\Lambda_1$. (iii) This is Lemma~\ref{lem:basic}(e).
\end{proof}

In the terminology of \citet{OgryczakRuszczynski2002}, \eqref{eq:areaPrimal} and \eqref{eq:areaDual} are the areas of the primal and the dual dispersion space of $L_t$. Thus $\xi$, $\Rrho$ and $\zeta_1$ are three norms of the same function $h$, which is concave in $x$ and vanishes on the boundary; in particular $\xi$ and $\zeta_1$, which are defined without rearrangement, are functionals of $\Cr$, compare \citet[Section~2]{SDS2024}. We have not found the identity $\Rrho=\Lambda_1$ in \citet{ALFT2026}, \citet{AnsariFuchs2025} or \citet{AnsariRockel2025}. It is consistent with the values $\Lambda_1=\alpha/3$ for the Farlie--Gumbel--Morgenstern family and $\Lambda_1=\alpha$ for the Fr\'echet family reported in \citet{ALFT2026}, for $\alpha\ge0$; for the former, $h=\alpha\,u(1-u)\,t(1-t)$ gives $12\int\!\!\int h=\alpha/3$.

\begin{corollary}\label{cor:ineq}
For every copula $C$,
\[
\tfrac23\,\zeta_1(C)^2\;\le\;\xi(C)\;\le\;\zeta_1(C)\;\le\;\Rrho(C)\;\le\;\min\bigl\{2\zeta_1(C)-\zeta_1(C)^2,\ \sqrt{2\,\xi(C)}\bigr\}.
\]
The constants $\tfrac23$ and $\sqrt2$ cannot be improved, and the set of pairs $(\zeta_1(C),\Rrho(C))$ is exactly $\{(\zeta,R):0\le\zeta\le1,\ \zeta\le R\le2\zeta-\zeta^2\}$. Moreover, $\xi(C)=\zeta_1(C)$ if and only if $L_t\in\{0,t,1\}$ almost surely for almost every $t$, and $\xi(C)=\Rrho(C)$ if and only if $C=\Pi$ or $C$ is completely dependent.
\end{corollary}

\begin{proof}
Let $m_t:=\max_xh_t(x)=\E(L_t-t)_+=\E(t-L_t)_+$, see Lemma~\ref{lem:basic}(e). For $\ell\in[0,1]$ we have $(\ell-t)^2\le(1-t)(\ell-t)_++t(t-\ell)_+$, with equality if and only if $\ell\in\{0,t,1\}$; hence $\Var(L_t)\le m_t$, and Theorem~\ref{thm:dictionary}(i),(iii) gives $\xi\le\zeta_1$ together with the equality case. Conversely, $m_t=\tfrac12\E|L_t-t|\le\tfrac12\Var(L_t)^{1/2}$, and Jensen's inequality gives $\zeta_1\le3\int_0^1\Var(L_t)^{1/2}\,\dd t\le3(\xi/6)^{1/2}$, that is, $\xi\ge\tfrac23\zeta_1^2$. A concave function $w$ on $[0,1]$ with $w(0)=w(1)=0$ and maximum $m$ satisfies $\int_0^1w\ge m/2$, which gives $\zeta_1\le\Rrho$. By Lemma~\ref{lem:basic}(b), $h_t\le\min\{x(1-t),m_t,(1-x)t\}$ and $m_t\le t(1-t)$, so that $\int_0^1h_t(x)\,\dd x\le m_t-m_t^2/(2t(1-t))$. By the Cauchy--Schwarz inequality, $(\int_0^1m_t\,\dd t)^2\le\int_0^1t(1-t)\,\dd t\int_0^1m_t^2/(t(1-t))\,\dd t$, hence $\Rrho\le12\int_0^1m_t\,\dd t-6\int_0^1m_t^2/(t(1-t))\,\dd t\le2\zeta_1-\zeta_1^2$. Integration by parts and the Cauchy--Schwarz inequality yield
\[
\int_0^1h_t(x)\,\dd x=-\int_0^1\bigl(x-\tfrac12\bigr)\,\partial_xh_t(x)\,\dd x\le\tfrac{1}{\sqrt{12}}\sqrt{\Var(L_t)},
\]
so that, by Jensen's inequality, $\Rrho\le\sqrt{12}\int_0^1\sqrt{\Var(L_t)}\,\dd t\le\sqrt{12\,\xi/6}=\sqrt{2\xi}$.

For the optimality of $\sqrt2$ let $k_\varepsilon(v)=\min\{v,1-v,\varepsilon\}$ with $\varepsilon\in(0,\tfrac12]$ and $C_\varepsilon(u,v)=uv+u(1-u)k_\varepsilon(v)$, which is a copula because $|k_\varepsilon'|\le1$. Then $h_t(x)=x(1-x)k_\varepsilon(t)$, hence $\Rrho(C_\varepsilon)=2\int k_\varepsilon=2\varepsilon(1-\varepsilon)$ and $\xi(C_\varepsilon)=2\int k_\varepsilon^2=2\varepsilon^2(1-\tfrac43\varepsilon)$, and $\Rrho/\sqrt{\xi}\to\sqrt2$ as $\varepsilon\to0$. For the optimality of $\tfrac23$ let $D_\varepsilon(u,v)=uv+\min\{u,1-u\}\,k_\varepsilon(v)$. Since $\partial_1D_\varepsilon(u,v)=v\pm k_\varepsilon(v)$ is nondecreasing in $v$ with values in $[0,1]$, $D_\varepsilon$ is a copula, and $L_t=t\pm k_\varepsilon(t)$ with probability $\tfrac12$ each. Hence $\xi(D_\varepsilon)=6\int k_\varepsilon^2=6\varepsilon^2(1-\tfrac43\varepsilon)$ and $\zeta_1(D_\varepsilon)=3\int k_\varepsilon=3\varepsilon(1-\varepsilon)$, so that $\xi(D_\varepsilon)/\zeta_1(D_\varepsilon)^2\to\tfrac23$.

For the region let $c\in[0,1]$, and let $P_c$ be given by $\partial_1P_c(u,t)=1$ for $u\le ct$, $\partial_1P_c(u,t)=t$ for $ct<u<1-c(1-t)$, and $\partial_1P_c(u,t)=0$ for $u\ge1-c(1-t)$. This function is nonincreasing in $u$, nondecreasing in $t$ and has mean $t$, so $P_c$ is a stochastically increasing copula; it is the rearrangement of the copula under which $V=U/c$ on $\{U<c\}$ and the conditional law of $V$ given $U=u$ is uniform for $u\ge c$. We have $h_t(x)=\min\{x(1-t),ct(1-t),(1-x)t\}$, hence $\zeta_1(P_c)=c$ and $\Rrho(P_c)=2c-c^2$, while the Fr\'echet copula $F_c:=(1-c)\Pi+cM$ satisfies $\zeta_1(F_c)=\Rrho(F_c)=c$. Both copulas are stochastically increasing, so that $h$ is affine along the segment $(1-s)F_c+sP_c$, $s\in[0,1]$, and for both $h_t$ attains its maximum at $x=t$. Hence $\zeta_1=c$ along the segment, while $\Rrho$ runs through $[c,2c-c^2]$.

If $\xi(C)=\Rrho(C)$, then $\xi=\zeta_1=\Rrho$; equality in $\zeta_1\le\Rrho$ forces $h_t$ to be a triangle, that is, $L_t$ takes at most two values, and together with $L_t\in\{0,t,1\}$ and $\E L_t=t$ this gives $\Var(L_t)\in\{0,t(1-t)\}$ for almost every $t$. Since $t\mapsto L_t$ is continuous in $L^1$ and bounded, $t\mapsto\Var(L_t)$ is continuous on $(0,1)$, so one of the two alternatives holds for all $t$.
\end{proof}

\begin{remark}[The exact upper bound]\label{rem:region}
For stochastically increasing copulas $\Cr=C$ and $\Rrho=\rho$, so the corollary contains the inequality $\xi\le\rho$ of \citet[Theorem~2]{AnsariRockel2025}, including the equality case. These authors also determine the exact region of the pairs $(\xi(C),\rho(C))$: it is given by $|\rho|\le M_\xi$ with an explicit function $M$, where $M_\xi=\beta-\tfrac3{10}\beta^2$ and $\xi=\tfrac1{10}\beta^2(5-2\beta)$ for a parameter $\beta\in(0,1]$, that is, for $\xi\le\tfrac3{10}$ \citep[Theorem~1]{AnsariRockel2025}. Since $\xi(\Cr)=\xi(C)$ and $\rho(C)\le\rho(\Cr)=\Rrho(C)$, their result applied to $\Cr$ shows that
\[
\Rrho(C)\le M_{\xi(C)}\qquad\text{for every copula }C,
\]
and that this bound is attained for every value of $\xi$, namely by the copulas on the right boundary of their region. Thus $M_\xi$ is the exact upper bound for $\Rrho$ in terms of $\xi$; Theorem~\ref{thm:exact} below gives a short direct proof and the corresponding bound for every spectrum. In the parametrisation above, $2\xi-M_\xi^2=\beta^3(20-9\beta)/100$, so that $M_\xi<\sqrt{2\xi}$ and $M_\xi/\sqrt{2\xi}\to1$ as $\xi\to0$; this is why the constant $\sqrt2$ in Corollary~\ref{cor:ineq} is optimal but not attained. The inequality $\xi\le\Lambda_1$ can also be read off from the monotonicity of $p\mapsto\Lambda_p$ in \citet{ALFT2026}. Equality $\zeta_1=\Rrho$ holds if and only if $L_t$ takes at most two values for almost every $t$.
\end{remark}

\section{Spectral dependence measures}\label{sec:spectral}

\subsection{One spectrum, two directions}

A \emph{spectrum} is a bounded, nonincreasing probability density $\phi$ on $(0,1)$ that is not constant; we write $\Fp(s)=\int_0^s\phi$ for its distribution function, which is concave with $\Fp(s)>s$ on $(0,1)$, and $\nu_\phi$ for the measure $-\dd\phi$ on $(0,1)$, whose total mass $\bar\nu_\phi:=\phi(0+)-\phi(1-)$ is positive and finite. For a random variable $X$ with values in $[0,1]$ let
\[
\mathcal M^-_\phi(X):=\int_0^1q_X(s)\,\phi(s)\,\dd s,\qquad \mathcal M^+_\phi(X):=\int_0^1X^*(s)\,\phi(s)\,\dd s .
\]
$\mathcal M^+_\phi$ is the spectral risk measure of \citet{Acerbi2002} if $X$ is a loss, and $\mathcal M^-_\phi$ is the corresponding valuation if $X$ is a payoff; both overweight unfavourable outcomes, and $\mathcal M^-_\phi(X)\le\E X\le \mathcal M^+_\phi(X)$. For $\phi_1(s)=2(1-s)$ one has $\mathcal M^-_{\phi_1}(X)=\E\min\{X,X'\}$ and $\mathcal M^+_{\phi_1}(X)=\E\max\{X,X'\}$ with an independent copy $X'$.

\begin{definition}\label{def:spectral}
Let $\phi$ be a spectrum, $N_\phi:=\int_0^1(\Fp(t)-t)\,\dd t$, and $(U,V)\sim C$. The \emph{spectral dependence measures} of $C$ with spectrum $\phi$ are
\[
\omega_\phi(C):=\frac{\E\bigl[\mathcal M^-_\phi(V\mid U)\bigr]-\mathcal M^-_\phi(V)}{N_\phi},\qquad
\kappa_\phi(C):=\frac{1}{N_\phi}\int_0^1\bigl(\mathcal M^+_\phi(L_t)-t\bigr)\,\dd t,
\]
where $\mathcal M^-_\phi(V\mid U)$ is $\mathcal M^-_\phi$ applied to the conditional law of $V$ given $U$.
\end{definition}

Thus $\omega_\phi$ measures how much the spectral value of $V$ increases on average when $U$ is known, and $\kappa_\phi$ measures the average spectral risk of the conditional probabilities $L_t$ in excess of their means. Both measures are unchanged if $\phi$ is replaced by a spectrum $(1-\vartheta)+\vartheta\phi$ with $\vartheta>0$, because $\Fp(t)-t$, the two numerators and $N_\phi$ are all multiplied by $\vartheta$ (for $\omega_\phi$ use $\E L_t=t$ and Proposition~\ref{prop:two}). For instance, every affine spectrum $s\mapsto1+\vartheta(1-2s)=(1-\vartheta)+\vartheta\phi_1(s)$, $\vartheta\in(0,1]$, defines the same measures as $\phi_1$.

\begin{proposition}[Representations]\label{prop:two}
For every spectrum $\phi$ and every copula $C$,
\begin{align}
\omega_\phi(C)&=\frac{1}{N_\phi}\int_0^1\bigl(\Fp(t)-\E\,\Fp(L_t)\bigr)\,\dd t=\frac{1}{N_\phi}\int_0^1\!\!\int_0^1\bigl(\Fp(t)-\Fp(\partial_1C(u,t))\bigr)\,\dd u\,\dd t,\label{eq:xiphi}\\
\kappa_\phi(C)&=\frac{1}{N_\phi}\int_0^1\!\!\int_{(0,1)}\bigl(\Cr-\Pi\bigr)(x,t)\,\nu_\phi(\dd x)\,\dd t .\label{eq:kappaphi}
\end{align}
In particular $\omega_{\phi_1}=\xi$ and $\kappa_{\phi_1}=\Rrho$.
\end{proposition}

\begin{proof}
If $S$ has distribution function $\Fp$ and $X$ has distribution function $F$, then $\mathcal M^-_\phi(X)=\E q_X(S)=\int_0^1\P(q_X(S)>v)\,\dd v=\int_0^1(1-\Fp(F(v)))\,\dd v$. Applying this to the conditional law of $V$ given $U$, whose distribution function at $v$ is $L_v$, and to the uniform law of $V$, and using Fubini's theorem, gives $\E[\mathcal M^-_\phi(V\mid U)]-\mathcal M^-_\phi(V)=\int_0^1(\Fp(v)-\E\Fp(L_v))\,\dd v$. For \eqref{eq:kappaphi} we integrate by parts: since $h_t$ is Lipschitz with $h_t(0)=h_t(1)=0$ and $\phi$ is bounded and of bounded variation, $\mathcal M^+_\phi(L_t)-t=\int_0^1\phi(s)\,\partial_sh_t(s)\,\dd s=\int_{(0,1)}h_t\,\dd\nu_\phi$. Finally $F_{\phi_1}(s)=2s-s^2$ gives $N_{\phi_1}=\tfrac16$ and $\omega_{\phi_1}(C)=6\int_0^1\Var(L_t)\,\dd t$, and $\nu_{\phi_1}$ is twice Lebesgue measure; now apply Theorem~\ref{thm:dictionary}.
\end{proof}

Thus $\xi$ is a spectral risk functional: $\tfrac16(\xi(C)+2)=\int\!\!\int(\partial_1C)^2=\E[\E(\min\{V,V'\}\mid U)]$ for conditionally independent copies $V,V'$, which is the average mean--Gini valuation of the conditional laws in the sense of \citet{Yitzhaki1982}; equivalently $\xi=1-3\,\E|V-V'|$, a representation also used by \citet{Limbach2026}. The two constructions are related by an exchange of axes: by \eqref{eq:xiphi}, a rank-dependent functional in the sense of \citet{Yaari1987} applied to the conditional laws of $V$ becomes an expected-utility functional $\E\,\Fp(L_t)$ of the conditional probabilities, because in the graph of $t\mapsto\partial_1C(u,t)$ outcomes and probabilities exchange their roles. Hence the measures $\omega_\phi$ are integral functionals of $\partial_1C$, whereas by \eqref{eq:kappaphi} the measures $\kappa_\phi$ are positive linear functionals of $\Cr-\Pi$, a subclass of the rearranged dependence measures of \citet{SDS2024}.

\begin{proposition}[Scalar products]\label{prop:pairing}
For a spectrum $\phi$ let $W_\phi(x,t):=\Fp(x)-x$. Then, for every copula $C$,
\begin{equation}\label{eq:pairing}
N_\phi\,\kappa_\phi(C)=\sprod{\Cr-\Pi}{W_\phi}_1=\int_0^1\!\!\int_0^1\bigl(\phi(x)-1\bigr)\,\partial_1(\Cr-\Pi)(x,t)\,\dd x\,\dd t,
\end{equation}
and consequently
\begin{equation}\label{eq:pairingCS}
\kappa_\phi(C)\le\frac{\norm{\phi-1}_{L^2}}{N_\phi}\,\sqrt{\frac{\xi(C)}{6}} .
\end{equation}
The constant in \eqref{eq:pairingCS} cannot be improved. For $\phi=\phi_1$ the inequality is $\Rrho\le\sqrt{2\xi}$.
\end{proposition}

\begin{proof}
We have $\partial_1W_\phi=\phi-1$ and $\int_0^1\partial_xh_t(x)\,\dd x=0$, so \eqref{eq:pairing} is the identity $\mathcal M^+_\phi(L_t)-t=\int_0^1\phi\,\partial_xh_t$ from the proof of Proposition~\ref{prop:two}, integrated over $t$. The Cauchy--Schwarz inequality and Theorem~\ref{thm:dictionary}(i) give \eqref{eq:pairingCS}, since $\norm{W_\phi}_1=\norm{\phi-1}_{L^2}$ and $\norm{\Cr-\Pi}_1^2=\xi(C)/6$. For $\phi_1$ we have $\norm{\phi_1-1}_{L^2}^2=\tfrac13$ and $N_{\phi_1}=\tfrac16$. For optimality let $m:=\norm{\phi-1}_\infty$, $\gamma:=(\phi-1)/m$, $\Gamma(u):=\int_0^u\gamma$, and, with $k_\varepsilon$ as in the proof of Corollary~\ref{cor:ineq}, $C_\varepsilon(u,v):=uv+\Gamma(u)k_\varepsilon(v)$. Since $|\gamma|\le1$, $\int_0^1\gamma=0$ and $\gamma$ is nonincreasing, this is a stochastically increasing copula with $\partial_xh_t(x)=\gamma(x)k_\varepsilon(t)$. Hence $N_\phi\kappa_\phi(C_\varepsilon)=m^{-1}\norm{\phi-1}_{L^2}^2\int k_\varepsilon$ and $\xi(C_\varepsilon)/6=m^{-2}\norm{\phi-1}_{L^2}^2\int k_\varepsilon^2$, so the quotient of the two sides of \eqref{eq:pairingCS} is $\int k_\varepsilon/(\int k_\varepsilon^2)^{1/2}=(1-\varepsilon)(1-\tfrac43\varepsilon)^{-1/2}\to1$ as $\varepsilon\to0$.
\end{proof}

Equality in \eqref{eq:pairingCS} would require $\partial_1(\Cr-\Pi)$ to be a multiple of $\phi-1$ that does not depend on $t$, which is impossible for a copula other than $\Pi$. For the tail spectra $\phi_p(s)=2(p-s)_+/p^2$ of Section~\ref{sec:sharp} one has $\norm{\phi_p-1}_{L^2}^2=\tfrac4{3p}-1$, so the control of $\kappa_{\phi_p}$ by $\xi$ degenerates like $p^{-1/2}$ as $p\to0$, and by the optimality statement this degeneration is genuine.

The Cauchy--Schwarz inequality ignores the constraints $0\le\partial_1\Cr\le1$. Taking them into account gives the exact upper bound for $\kappa_\phi$ in terms of $\xi$, for every spectrum. For $\lambda>0$ and $\mu\in\mathbb R$ let $\ell_{\lambda,\mu}(x):=\min\{1,\max\{0,\lambda(\phi(x)-1-\mu)\}\}$. The map $\mu\mapsto\int_0^1\ell_{\lambda,\mu}$ is continuous and nonincreasing with values $1$ for small and $0$ for large $\mu$, and the functions $\ell_{\lambda,\mu}$ are pointwise nonincreasing in $\mu$. Hence for every $t\in[0,1]$ there is a function
\begin{equation}\label{eq:extremal}
\ell_\lambda(\cdot,t)=\ell_{\lambda,\mu_\lambda(t)}\qquad\text{with}\quad\int_0^1\ell_\lambda(x,t)\,\dd x=t,
\end{equation}
it is unique almost everywhere, and $\ell_\lambda(\cdot,t)\le\ell_\lambda(\cdot,t')$ for $t<t'$; we choose $\mu_\lambda(t)$ as the largest admissible value in the interval $[\inf\phi-1-1/\lambda,\ \sup\phi-1]$, which is nonincreasing in $t$. Since moreover $\ell_\lambda(\cdot,0)=0$, $\ell_\lambda(\cdot,1)=1$ and $\ell_\lambda$ is nonincreasing in $x$,
\[
C_{\phi,\lambda}(x,t):=\int_0^x\ell_\lambda(s,t)\,\dd s
\]
is a stochastically increasing copula with $\partial_1C_{\phi,\lambda}=\ell_\lambda$. Let further $(a_i,c_i)$ be the maximal open subintervals of $(0,1)$ on which $\phi$ is constant, with lengths $w_i:=c_i-a_i$, and let $O_\phi$ be the ordinal sum of $\Pi$ on the squares $[a_i,c_i]^2$ and of $M$ elsewhere; explicitly, $\partial_1O_\phi(u,t)=\mathbf 1\{u\le t\}$ if $t$ lies in none of these intervals, and $\partial_1O_\phi(u,t)=\mathbf 1\{u\le a_i\}+\frac{t-a_i}{w_i}\mathbf 1\{a_i<u<c_i\}$ if $t\in(a_i,c_i)$. If $\phi$ is strictly decreasing, then $O_\phi=M$.

\begin{theorem}[Exact upper bound]\label{thm:exact}
Let $\phi$ be a spectrum, and for $\lambda>0$ let $K(\lambda):=N_\phi\,\kappa_\phi(C_{\phi,\lambda})$ and $X(\lambda):=\xi(C_{\phi,\lambda})$.
\begin{enumerate}
\item[\textup{(a)}] For every copula $C$ and every $\lambda>0$,
\begin{equation}\label{eq:lagrange}
N_\phi\,\kappa_\phi(C)-\frac{\xi(C)}{12\,\lambda}\;\le\;K(\lambda)-\frac{X(\lambda)}{12\,\lambda},
\end{equation}
with equality if and only if $\Cr=C_{\phi,\lambda}$. In particular $\kappa_\phi(C)\le\kappa_\phi(C_{\phi,\lambda})$ whenever $\xi(C)\le\xi(C_{\phi,\lambda})$.
\item[\textup{(b)}] $K$ and $X$ are continuous and nondecreasing, $X(\lambda)\to0$ as $\lambda\to0$, and $K(\lambda)\to N_\phi$ as $\lambda\to\infty$. Let $\xi_\phi^\infty:=\lim_{\lambda\to\infty}X(\lambda)$. The function
\[
M_\phi(\xi):=\sup\{\kappa_\phi(C):\ \xi(C)=\xi\},\qquad\xi\in(0,1],
\]
satisfies $M_\phi(X(\lambda))=\kappa_\phi(C_{\phi,\lambda})$ for every $\lambda>0$, is concave and nondecreasing on $(0,\xi_\phi^\infty)$, and equals $1$ on $[\xi_\phi^\infty,1]$. The supremum is attained for every $\xi$. For $\xi\in(0,\xi_\phi^\infty)$ it is attained exactly by the copulas $C$ with $\Cr=C_{\phi,\lambda}$, where $X(\lambda)=\xi$; this copula does not depend on the choice of $\lambda$.
\item[\textup{(c)}] $\xi_\phi^\infty=\min\{\xi(C):\kappa_\phi(C)=1\}=1-\sum_iw_i^2$, and the minimum is attained exactly by the copulas $C$ with $\Cr=O_\phi$. In particular $\xi_\phi^\infty=1$ if and only if $\phi$ is strictly decreasing.
\end{enumerate}
\end{theorem}

\begin{proof}
(a) Let $w:=\phi-1$. For fixed $(x,t)$ the function $y\mapsto w(x)y-y^2/(2\lambda)-\mu_\lambda(t)y$ is strictly concave, and its unique maximiser on $[0,1]$ is $\ell_\lambda(x,t)$. Let $C$ be a copula and $\ell:=\partial_1\Cr$. Then $0\le\ell\le1$, so that
\[
w(x)\ell(x,t)-\frac{\ell(x,t)^2}{2\lambda}-\mu_\lambda(t)\ell(x,t)\;\le\;w(x)\ell_\lambda(x,t)-\frac{\ell_\lambda(x,t)^2}{2\lambda}-\mu_\lambda(t)\ell_\lambda(x,t)
\]
for all $(x,t)$. We integrate over the unit square. Since $\int_0^1\ell(x,t)\,\dd x=t=\int_0^1\ell_\lambda(x,t)\,\dd x$, the terms containing $\mu_\lambda$ cancel. By \eqref{eq:pairing}, $\int\!\!\int w\ell=N_\phi\kappa_\phi(C)$ and $\int\!\!\int w\ell_\lambda=K(\lambda)$, and $\int\!\!\int\ell^2=\tfrac16(\xi(C)+2)$ by Theorem~\ref{thm:dictionary}(i), likewise for $\ell_\lambda$. This is \eqref{eq:lagrange}, and equality holds if and only if $\ell=\ell_\lambda$ almost everywhere.

(b) Applying \eqref{eq:lagrange} to $C=C_{\phi,\lambda'}$ gives $K(\lambda')-K(\lambda)\le(X(\lambda')-X(\lambda))/(12\lambda)$, and exchanging $\lambda$ and $\lambda'$ gives $K(\lambda')-K(\lambda)\ge(X(\lambda')-X(\lambda))/(12\lambda')$. Hence $(X(\lambda')-X(\lambda))(1/\lambda-1/\lambda')\ge0$, so $X$ and then $K$ are nondecreasing. Moreover $X(\lambda)=X(\lambda')$ implies $K(\lambda)=K(\lambda')$, so that $C_{\phi,\lambda'}$ attains equality in \eqref{eq:lagrange} for the parameter $\lambda$ and therefore coincides with $C_{\phi,\lambda}$. If $\lambda_n\to\lambda$, then for fixed $t$ the bounded sequence $\mu_{\lambda_n}(t)$ has limit points, and for each of them $\ell_{\lambda,\mu}$ has integral $t$ by dominated convergence, so that $\ell_{\lambda_n}(\cdot,t)\to\ell_\lambda(\cdot,t)$ in $L^1$; by dominated convergence in $t$, $X$ and $K$ are continuous. For $C=\Pi$ and $C=M$, \eqref{eq:lagrange} reads $X(\lambda)\le12\,\lambda\,K(\lambda)\le12\,\lambda\,N_\phi$ and $K(\lambda)\ge N_\phi-(1-X(\lambda))/(12\lambda)$, which gives the two limits. By \eqref{eq:lagrange}, $N_\phi\kappa_\phi(C)\le K(\lambda)$ if $\xi(C)=X(\lambda)$, with equality if and only if $\Cr=C_{\phi,\lambda}$, and by continuity every $\xi\in(0,\xi_\phi^\infty)$ is of the form $X(\lambda)$. Moreover \eqref{eq:lagrange} shows that $N_\phi M_\phi(\xi)=\inf_{\lambda>0}\{K(\lambda)+(\xi-X(\lambda))/(12\lambda)\}$ on this interval, an infimum of nondecreasing affine functions. Finally, by (c), $\xi(O_\phi)=\xi_\phi^\infty$ and $\kappa_\phi(O_\phi)=1$. The copulas $(1-\sigma)O_\phi+\sigma M$, $\sigma\in[0,1]$, are stochastically increasing, so $\kappa_\phi$ is affine along this segment by \eqref{eq:kappaphi} and equals $1$ at both ends, while $\xi$ is continuous and takes all values in $[\xi_\phi^\infty,1]$. Hence $M_\phi=1$ on $[\xi_\phi^\infty,1]$, with attainment.

(c) Let $\xi^*:=1-\sum_iw_i^2$, and let $C$ be a copula with $\kappa_\phi(C)=1$. Since $0\le h_t\le h^M_t$ with $h^M_t(x)=M(x,t)-xt$ and $\int_0^1\!\int_{(0,1)}h^M_t\,\dd\nu_\phi\,\dd t=N_\phi$ (see the proof of Theorem~\ref{thm:mirror}(a)), $\kappa_\phi(C)=1$ implies that for almost every $t$ the continuous function $h^M_t-h_t\ge0$ vanishes on the support of $\nu_\phi$, which is the complement in $(0,1)$ of the union of the intervals $(a_i,c_i)$. If $t$ lies in this support, then $h_t(t)=h^M_t(t)$, so $L_t\in\{0,1\}$ and $\Var(L_t)=t(1-t)$. If $t\in(a_i,c_i)$, then $h_t=h^M_t$ at $a_i$ and at $c_i$, that is, $L_t^*=1$ on $[0,a_i)$ and $L_t^*=0$ on $(c_i,1]$. Hence $\int_{a_i}^{c_i}L_t^*=t-a_i$, and the Cauchy--Schwarz inequality gives $\E L_t^2\ge a_i+(t-a_i)^2/w_i$, with equality if and only if $L_t^*=(t-a_i)/w_i$ on $(a_i,c_i)$; equivalently, $\Var(L_t)\ge t(1-t)-(t-a_i)(c_i-t)/w_i$. Since $\int_{a_i}^{c_i}(t-a_i)(c_i-t)\,\dd t/w_i=w_i^2/6$, integration gives $\xi(C)\ge\xi^*$, with equality if and only if every $L_t$ has the same law as under $O_\phi$, that is, if and only if $\Cr=O_\phi$. Under $O_\phi$ we have $h_t=h^M_t$ on $[0,a_i]\cup[c_i,1]$ for $t\in(a_i,c_i)$ and $L_t\in\{0,1\}$ for the other $t$, so $\kappa_\phi(O_\phi)=1$ by \eqref{eq:kappaphi}. Thus $\xi^*=\min\{\xi(C):\kappa_\phi(C)=1\}$. Next, \eqref{eq:lagrange} for a copula $C$ with $\kappa_\phi(C)=1$, together with $K(\lambda)\le N_\phi$, gives $X(\lambda)\le\xi(C)$ for all $\lambda$, hence $\xi_\phi^\infty\le\xi^*$. Conversely let $\lambda_n\to\infty$. Since the set of copulas is compact with respect to uniform convergence, a subsequence of $C_{\phi,\lambda_n}$ converges uniformly to a copula $D$, which is stochastically increasing, and $\kappa_\phi(D)=\lim K(\lambda_n)/N_\phi=1$ by \eqref{eq:kappaphi}. Moreover $\partial_1C_{\phi,\lambda_n}\to\partial_1D$ weakly in $L^2([0,1]^2)$, because the sequence is bounded by $1$ and $\int\!\!\int\partial_1C\,\mathbf 1_{(u_1,u_2]\times(t_1,t_2]}=\int_{t_1}^{t_2}(C(u_2,t)-C(u_1,t))\,\dd t$ converges for every rectangle. By weak lower semicontinuity of the norm, $\xi(D)\le\liminf X(\lambda_n)=\xi_\phi^\infty$, hence $\xi^*\le\xi_\phi^\infty$. Finally, $\sum_iw_i^2=0$ if and only if $\phi$ is constant on no interval, that is, strictly decreasing.
\end{proof}

By \eqref{eq:lagrange}, $C_{\phi,\lambda}$ maximises $N_\phi\kappa_\phi-\xi/(12\lambda)$ over all copulas; the proof uses only the constraints $0\le\partial_1\Cr\le1$ and $\int_0^1\partial_1\Cr(x,t)\,\dd x=t$, and $1/(12\lambda N_\phi)$ is the slope of the concave boundary curve $M_\phi$ at $\xi=X(\lambda)$. Near the origin, $M_\phi(\xi)\le N_\phi^{-1}\norm{\phi-1}_{L^2}\sqrt{\xi/6}$ by Proposition~\ref{prop:pairing}, with optimal constant. By part (c), if $\phi$ is not strictly decreasing, $\kappa_\phi$ reaches its maximal value already at the level $\xi_\phi^\infty<1$; compare Theorem~\ref{thm:mirror}(d).

\begin{example}[Two explicit boundary curves]\label{ex:exact}
\leavevmode\par\noindent (i) For $\phi=\phi_1$ we have $\ell_\lambda(x,t)=\min\bigl\{1,\max\{0,\beta(s(t)-x)\}\bigr\}$ with $\beta=2\lambda$, so that $C_{\phi_1,\lambda}$ has a diagonal band structure. A direct computation gives
\[
\bigl(\xi,\Rrho\bigr)(C_{\phi_1,\lambda})=
\begin{cases}
\bigl(\tfrac1{10}\beta^2(5-2\beta),\ \beta-\tfrac3{10}\beta^2\bigr),&\beta\le1,\\[0.5ex]
\bigl(1-\tfrac1\beta+\tfrac3{10\beta^2},\ 1-\tfrac1{2\beta^2}+\tfrac1{5\beta^3}\bigr),&\beta\ge1.
\end{cases}
\]
This is the boundary curve $M_\xi$ of the exact $(\xi,\rho)$ region of \citet[Theorem~1]{AnsariRockel2025}, see Remark~\ref{rem:region}, and by their uniqueness statement the copulas $C_{\phi_1,\lambda}$ are their boundary copulas. Since $\rho(C)\le\Rrho(C)$, Theorem~\ref{thm:exact} contains the upper boundary of that region, with a short proof.

(ii) For the spectrum $\phi=p^{-1}\mathbf 1_{[0,p)}$ of expected shortfall, $\ell_\lambda(\cdot,t)$ takes two values, and one finds $C_{\phi,\lambda}=\Pi+A_p\otimes g_\lambda$ with $A_p$ as in Section~\ref{sec:sharp} and the plateau function $g_\lambda(t)=\min\{(1-p)\lambda,(1-p)t,p(1-t)\}$. For $\lambda\ge p$ this is the block copula $B_p$ of Section~\ref{sec:stoploss}, which is $O_\phi$. With $\sigma:=\min\{\lambda/p,1\}$,
\[
\kappa_\phi(C_{\phi,\lambda})=\sigma(2-\sigma),\qquad\xi(C_{\phi,\lambda})=2p(1-p)\,\sigma^2(3-2\sigma),
\]
so that $\xi_\phi^\infty=2p(1-p)=1-p^2-(1-p)^2$, in accordance with Theorem~\ref{thm:exact}(c), and $M_\phi(\xi)=\sigma(2-\sigma)$ where $\sigma\in[0,1]$ solves $\sigma^2(3-2\sigma)=\xi/(2p(1-p))$. 
\end{example}

Since the measures $\kappa_\phi$ are linear functionals of $\Cr$, the same pointwise idea, without the quadratic term, determines the joint range of two of them.

\begin{proposition}[Two spectral measures]\label{prop:pairs}
Let $w$ be a bounded measurable function on $(0,1)$, $W(x):=\int_0^xw$, and let $\widehat W$ be the least concave majorant of $W$ on $[0,1]$. Then
\begin{equation}\label{eq:hull}
\sup_C\int_0^1\!\!\int_0^1w(x)\,\partial_1\Cr(x,t)\,\dd x\,\dd t=\int_0^1\widehat W(t)\,\dd t,
\end{equation}
and the supremum is attained by the ordinal sum of $\Pi$ on the squares $[a,c]^2$, where $(a,c)$ runs through the components of $\{\widehat W>W\}$, and of $M$ elsewhere. Consequently, for two spectra $\phi$ and $\psi$, the set $\mathcal R_{\phi,\psi}:=\{(\kappa_\phi(C),\kappa_\psi(C)):\ C\text{ a copula}\}$ is convex and compact, and for all $\alpha,\beta\in\mathbb R$
\[
\max_{(k,k')\in\mathcal R_{\phi,\psi}}\bigl(\alpha N_\phi k+\beta N_\psi k'\bigr)=\int_0^1\widehat W_{\alpha,\beta}(t)\,\dd t,\qquad W_{\alpha,\beta}:=\alpha(\Fp-\mathrm{id})+\beta(F_\psi-\mathrm{id}),
\]
where the maximum is attained by the ordinal sum associated with $W_{\alpha,\beta}$ as above.
\end{proposition}

\begin{proof}
For fixed $t$ the function $\ell:=\partial_1\Cr(\cdot,t)$ is nonincreasing with values in $[0,1]$ and $\int_0^1\ell=t$, so there is a random variable $S$ with values in $[0,1]$ such that $\ell(x)=\P(S>x)$ for almost every $x$, and $\E S=t$. By Fubini's theorem $\int_0^1w\ell=\E W(S)\le\E\widehat W(S)\le\widehat W(t)$ by Jensen's inequality. This gives the inequality in \eqref{eq:hull}. Since $W$ is Lipschitz, $\widehat W$ is continuous and coincides with $W$ at $0$ and $1$, the contact set $Z:=\{\widehat W=W\}$ is closed, and $\widehat W$ is affine on the closure of each component $(a,c)$ of $[0,1]\setminus Z$. Let $O$ be the ordinal sum in the statement. For $t\in Z$, $\partial_1O(\cdot,t)=\mathbf 1_{[0,t]}$ corresponds to $S\equiv t$, and $\E W(S)=W(t)=\widehat W(t)$. For $t\in(a,c)$, $\partial_1O(\cdot,t)$ corresponds to a variable $S$ with values in $\{a,c\}$ and mean $t$, and $\E W(S)=\widehat W(t)$ because $\widehat W$ is affine on $[a,c]$ and coincides with $W$ at $a$ and $c$. Hence $O$ attains the supremum. For the second statement recall that $\kappa_\phi$ and $\kappa_\psi$ depend on $C$ only through $\Cr$, which runs through the set of all stochastically increasing copulas. This set is convex and compact with respect to uniform convergence, and by \eqref{eq:kappaphi} both measures are affine and continuous on it, so $\mathcal R_{\phi,\psi}$ is convex and compact. By \eqref{eq:pairing}, and since $\int_0^1(\phi-1)=0$, the functional $\alpha N_\phi\kappa_\phi+\beta N_\psi\kappa_\psi$ is the left-hand side of \eqref{eq:hull} with $w=\alpha(\phi-1)+\beta(\psi-1)$, whose integral is $W_{\alpha,\beta}$.
\end{proof}

\begin{example}\label{ex:pairs}
Let $\phi=p^{-1}\mathbf 1_{[0,p)}$ be the spectrum of expected shortfall and $\psi=\phi_1$, and take $\alpha=1$, $\beta=-1/p$. Then $W_{\alpha,\beta}(x)=x(x-p)/p$ for $x\le p$ and $W_{\alpha,\beta}(x)=(1-x)(1-x/p)$ for $x\ge p$, so that $W_{\alpha,\beta}\le0$ with equality at $0$, $p$ and $1$, and $\widehat W_{\alpha,\beta}=0$. Since $N_\phi=\tfrac12(1-p)$ and $N_{\phi_1}=\tfrac16$, Proposition~\ref{prop:pairs} gives
\[
\kappa_\phi(C)\le\frac{\Rrho(C)}{3p(1-p)}\qquad\text{for every copula }C,
\]
with equality for $\Pi$, for the block copula $B_p$ and for their mixtures.
\end{example}

\begin{remark}[Relation to known classes]\label{rem:lit}
(a) Since $\E L_t=t$, the integrand in \eqref{eq:xiphi} can be written as $\E\,B_\phi(L_t,t)$ with the Bregman divergence $B_\phi(\ell,t)=\Fp(t)-\Fp(\ell)+\phi(t)(\ell-t)\ge0$ of the convex function $-\Fp$. The measures $\xi_\Psi$ of \citet{AnsariFuchs2025} and \citet{AnsariWiesel2026} are instead generated by $\Psi(\ell-t)$ with a convex function $\Psi$ of the difference, and the measures $\Lambda_\Psi$ of \citet{ALFT2026} by $\Psi(\ell-\ell')$. A Bregman divergence that depends only on the difference $\ell-t$ is generated by a quadratic function, so the integrands coincide only for the affine spectra, for which $B_\phi(\ell,t)=\vartheta(\ell-t)^2$; all of them give Chatterjee's $\xi$. For the difference-type measures, strict convexity of $\Psi$ at $0$ yields both characterisations at once \citep[Theorem~3.2]{ALFT2026} \citep[Theorem~4.1]{AnsariFuchs2025}; Theorem~\ref{thm:mirror} below shows that the situation for $\omega_\phi$ is different.

(b) If $(U^\uparrow,V^\uparrow)\sim\Cr$, then $\int_0^1\mathcal M^+_\phi(L_t)\,\dd t=\E[\phi(U^\uparrow)(1-V^\uparrow)]$, so that $N_\phi\,\kappa_\phi(C)=\int \chi_\phi\,\dd\Cr-\int \chi_\phi\,\dd\Pi$ with the supermodular function $\chi_\phi(u,v)=\phi(u)(1-v)$. Standardised linear functionals of this kind are mentioned as examples of rearranged dependence measures $R_\eta(C)=\eta(\Cr)$ in \citet{AnsariWiesel2026}, where the characterisations of independence and of perfect dependence are properties required of $\eta$ on stochastically increasing copulas. Theorem~\ref{thm:mirror} determines for which spectra they hold. The counterexample in part (d) of that theorem is of the same kind as the ordinal sum by which \citet{SDS2024} show that Blomqvist's $\beta$ does not characterise complete dependence on stochastically increasing copulas.

(c) The measures $\omega_\phi$ fit into the general framework of \citet{BorgonovoEtAl2025}, who study functionals $\E[\zeta(\P^V,\P^{V\mid U})]$ for separation measurements $\zeta$ that are convex in the second argument: take $\zeta(\nu,\nu')=\mathcal M^-_\phi(\nu')-\mathcal M^-_\phi(\nu)$, which is convex in $\nu'$ by the formula in the proof of Proposition~\ref{prop:two} (we do not check the lower semicontinuity that is part of their Definition~2.1, since we give direct proofs). They show that such a functional vanishes only under independence if and only if $\zeta$ is strict, and that it is maximal only under functional dependence if and only if $\zeta$ is strictly convex on Dirac masses \citep[Theorems~2.9 and 2.12]{BorgonovoEtAl2025}. For $\omega_\phi$ the second condition reads $\Fp(\tfrac12)>\tfrac12$ and holds for every spectrum, and by Theorem~\ref{thm:mirror} the first holds exactly for the strictly decreasing spectra. The statements on $\omega_\phi$ in Theorem~\ref{thm:mirror}(b)--(d) are therefore instances of their general principle, for which we give direct proofs. The measures $\kappa_\phi$ are not of this form, since they are not averages of a functional of the conditional laws of $V$. What we believe to be new is the common parametrisation of the two classes by a spectrum and the resulting symmetry of the two sets of conditions.
\end{remark}

\subsection{Mirror properties}\label{sec:mirror}

\begin{theorem}\label{thm:mirror}
Let $\phi$ be a spectrum.
\begin{enumerate}
\item[\textup{(a)}] $\omega_\phi$ and $\kappa_\phi$ take values in $[0,1]$, are monotone in the Schur order, depend on $C$ only through $\Cr$, and satisfy the information inequality: if $(U,V)\sim C$ and $f(U)$ has a continuous distribution, then the copula $C_f$ of $(f(U),V)$ satisfies $\omega_\phi(C_f)\le\omega_\phi(C)$ and $\kappa_\phi(C_f)\le\kappa_\phi(C)$. Moreover, both are convex functionals of $C$.
\item[\textup{(b)}] $\omega_\phi(C)=1$ if and only if $C$ is completely dependent, and $\kappa_\phi(C)=0$ if and only if $C=\Pi$.
\item[\textup{(c)}] If $\phi$ is strictly decreasing, then $\omega_\phi(C)=0$ if and only if $C=\Pi$, and $\kappa_\phi(C)=1$ if and only if $C$ is completely dependent.
\item[\textup{(d)}] If $\phi$ is not strictly decreasing, then both statements in \textup{(c)} fail: there are copulas $C\ne\Pi$ with $\omega_\phi(C)=0$ and copulas $C$ that are not completely dependent with $\kappa_\phi(C)=1$.
\end{enumerate}
\end{theorem}

\begin{proof}
(a) Since $\Fp$ is concave and $\Fp(\ell)\ge\ell$ on $[0,1]$, Jensen's inequality and $\E L_t=t$ give $0\le\Fp(t)-\E\Fp(L_t)\le\Fp(t)-t$. By Lemma~\ref{lem:basic}(b), $0\le h_t\le h^M_t$ with $h^M_t(x)=M(x,t)-xt$, and $\int_{(0,1)}h^M_t\,\dd\nu_\phi=\mathcal M^+_\phi(I_t)-t=\Fp(t)-t$ for a Bernoulli variable $I_t$ with mean $t$. This gives the ranges. Both measures depend only on the laws of the variables $L_t$, that is, on $\Cr$. Monotonicity follows from \eqref{eq:schur}, since $\ell\mapsto-\Fp(\ell)$ is convex and $h$ is pointwise monotone. If $U'=f(U)$, then $\P(V\le t\mid U')=\E[L_t\mid U']\le_{\cx}L_t$. Finally, $C\mapsto\partial_1C$ is linear, $-\Fp$ is convex, and $S_t(x)$ is, by \eqref{eq:supform}, a supremum of linear functionals of $\partial_1C(\cdot,t)$.

(b) $\omega_\phi(C)=1$ means $\E\Fp(L_t)=t=\E L_t$ for almost every $t$. Since $\Fp(\ell)>\ell$ on $(0,1)$, this holds if and only if $L_t\in\{0,1\}$ almost surely. If $\kappa_\phi(C)=0$, then $h_t=0$ $\nu_\phi$-almost everywhere for almost every $t$; a nonnegative concave function vanishing at $0$, at $1$ and at an interior point vanishes identically, so $L_t$ is constant. In both cases the property extends from almost every $t$ to all $t$ by continuity of $t\mapsto\Var(L_t)$.

(c) If $\phi$ is strictly decreasing, then $\Fp$ is strictly concave, so equality in Jensen's inequality forces $L_t=t$ almost surely. Moreover the support of $\nu_\phi$ is $[0,1]$, so $\kappa_\phi(C)=1$ forces the continuous functions $h_t$ and $h^M_t$ to coincide for almost every $t$; at $x=t$ this gives $S_t(t)=t=\E L_t$, hence $L_t^*=1$ on $[0,t)$ and $L_t\in\{0,1\}$. In both cases the property holds for almost every $t$ and extends to all $t$ by continuity of $t\mapsto\Var(L_t)$, as in (b).

(d) A nonincreasing function that is not strictly decreasing is constant on some interval $(a,c)$ with $0<a<c<1$. Then $\Fp$ is affine on $[a,c]$. Let $k(t)=\tfrac12\min\{(t-a)_+,(c-t)_+\}$ and $C(u,v)=uv+u(1-u)k(v)$, a copula different from $\Pi$ with $L_t=t+(1-2U)k(t)$. For $t\in(a,c)$ the variable $L_t$ takes values in $[a,c]$, and for $t\notin(a,c)$ it is constant; in both cases $\E\Fp(L_t)=\Fp(t)$, so $\omega_\phi(C)=0$. Next let $C$ be the ordinal sum of $M$, $\Pi$, $M$ with respect to the partition $[0,a],[a,c],[c,1]$. For $t\notin(a,c)$ we have $L_t\in\{0,1\}$, and for $t\in(a,c)$ we have $L_t=1$ on $\{U<a\}$, $L_t=0$ on $\{U>c\}$ and $L_t=(t-a)/(c-a)$ otherwise. Hence $h_t=h^M_t$ on $[0,a]\cup[c,1]$ for all $t$. Since $\nu_\phi((a,c))=0$, we get $\kappa_\phi(C)=1$, although $C$ is not completely dependent.
\end{proof}

For the spectrum $\phi=p^{-1}\mathbf 1_{[0,p)}$ of expected shortfall, $\omega_\phi(C)=0$ holds if and only if $L_p=p$ almost surely, that is, if and only if the single event $\{V\le p\}$ is independent of $U$, and $\kappa_\phi(C)$ is a normalised integral of $t\mapsto\Cr(p,t)-pt$.

\begin{proposition}[Continuity]\label{prop:cont}
For every spectrum $\phi$ and all copulas $C,D$,
\[
|\omega_\phi(C)-\omega_\phi(D)|\le\frac{\bar\nu_\phi}{2N_\phi}\,D_1(C,D),\qquad|\kappa_\phi(C)-\kappa_\phi(D)|\le\frac{\bar\nu_\phi}{2N_\phi}\,D_1(C,D).
\]
In particular $\max\{\omega_\phi(C),\kappa_\phi(C)\}\le\bar\nu_\phi\,\zeta_1(C)/(6N_\phi)$.
\end{proposition}

\begin{proof}
Realise $L^C_t$ and $L^D_t$ on the same variable $U$ and let $c_\phi:=\tfrac12(\phi(0+)+\phi(1-))$, so that $|\phi-c_\phi|\le\bar\nu_\phi/2$. Since $\E L^C_t=\E L^D_t$, we have $\E\Fp(L^C_t)-\E\Fp(L^D_t)=\E[\tilde F_\phi(L^C_t)-\tilde F_\phi(L^D_t)]$ with $\tilde F_\phi(\ell):=\Fp(\ell)-c_\phi\ell$, which is Lipschitz with constant $\bar\nu_\phi/2$; with \eqref{eq:xiphi} this gives the first inequality. Similarly, $\int_0^1(L^{C*}_t-L^{D*}_t)=0$ gives $\mathcal M^+_\phi(L^C_t)-\mathcal M^+_\phi(L^D_t)=\int_0^1(\phi-c_\phi)(L^{C*}_t-L^{D*}_t)$, and decreasing rearrangement is a contraction in $L^1$, so that $|\mathcal M^+_\phi(L^C_t)-\mathcal M^+_\phi(L^D_t)|\le\tfrac12\bar\nu_\phi\,\E|L^C_t-L^D_t|$. Integration over $t$ gives the second inequality, and the last statement is the case $D=\Pi$ together with $\zeta_1=3D_1(\cdot,\Pi)$.
\end{proof}

For $\phi=\phi_1$ the second inequality reads $|\Rrho(C)-\Rrho(D)|\le6D_1(C,D)$.

\subsection{Variability of the conditional laws}

The Schur order compares, for fixed $t$, the variability of $u\mapsto\partial_1C(u,t)$. Since $\omega_\phi$ is defined through the conditional laws $K_C(u,\cdot)$ of $V$, it is also monotone with respect to the variability of these laws for fixed $u$.

\begin{proposition}\label{prop:ssd}
Let $C$ and $D$ be copulas such that $K_C(u,\cdot)\le_{\cx}K_D(u,\cdot)$ for almost every $u$. Then $\omega_\phi(C)\ge\omega_\phi(D)$ for every spectrum $\phi$; in particular $\xi(C)\ge\xi(D)$. The assumption holds if and only if, for almost every $u$, the law $K_C(u,\cdot)$ dominates $K_D(u,\cdot)$ in second-order stochastic dominance.
\end{proposition}

\begin{proof}
For a probability measure $\mu$ on $[0,1]$ with quantile function $q_\mu$ and mean $m_\mu$ let $Q_\mu(p):=\int_0^pq_\mu$. Integration by parts gives
\begin{equation}\label{eq:Mminus}
\mathcal M^-_\phi(\mu)=\phi(1-)\,m_\mu+\int_{(0,1)}Q_\mu(p)\,\nu_\phi(\dd p).
\end{equation}
Second-order stochastic dominance of $\mu$ over $\mu'$ is equivalent to $Q_\mu\ge Q_{\mu'}$ on $[0,1]$ \citep[Theorem~3.2]{OgryczakRuszczynski2002}; in particular $m_\mu\ge m_{\mu'}$, and $\mu\le_{\cx}\mu'$ holds if and only if $\mu$ dominates $\mu'$ and the means are equal. Hence, under the assumption, \eqref{eq:Mminus} gives $\mathcal M^-_\phi(K_C(u,\cdot))\ge\mathcal M^-_\phi(K_D(u,\cdot))$ for almost every $u$, and integration over $u$ gives the claim, because the second term in the numerator of $\omega_\phi$ does not depend on the copula. For the last statement note that $\int_0^1m_{K_C(u,\cdot)}\,\dd u=\E V=\tfrac12$ for every copula; if the conditional means under $C$ are almost everywhere at least as large as those under $D$, they therefore coincide almost everywhere.
\end{proof}

For $\phi=\phi_1$ the proposition says that, among the distribution functions $F$ on $[0,1]$ with a given mean, $\int_0^1F(v)^2\,\dd v$ is nonincreasing in the convex order. The order in Proposition~\ref{prop:ssd} is restrictive: two comparable copulas have the same regression function $u\mapsto\E(V\mid U=u)$, so that, for instance, a copula is comparable with $\Pi$ only if this function is constant. Orders of this type are discussed in \citet[Section~5.2]{AnsariFuchs2025}. Whereas $\omega_\phi$ increases with the variability of the functions $u\mapsto\partial_1C(u,t)$, it decreases with the variability of the laws with distribution functions $t\mapsto\partial_1C(u,t)$.

\subsection{Behaviour near independence}

To compare the two classes we use a family of copulas in which all $L_t$ are affine images of a single random variable. Let $G$ be a random variable with $\E G=0$ and $|G|\le1$ that is not constant, let $\gamma$ be its decreasing quantile function, $\Gamma_G(u)=\int_0^u\gamma$, and for $\varepsilon\in[0,1]$
\begin{equation}\label{eq:lab}
C_{\varepsilon G}(u,v):=uv+\varepsilon\,\Gamma_G(u)\,v(1-v).
\end{equation}
This is a stochastically increasing copula of the type studied by \citet{RodriguezUbeda2004}, with $L_t=t+\varepsilon\,t(1-t)\,\gamma(U)$.

\begin{proposition}\label{prop:order}
Let $\phi,\psi$ be spectra and $G$ as above. Then, for $\varepsilon\in[0,1]$,
\[
\kappa_\psi(C_{\varepsilon G})=\frac{\mathcal M^+_\psi(G)}{6N_\psi}\,\varepsilon\quad\text{with }\mathcal M^+_\psi(G)>0,\qquad\omega_\phi(C_{\varepsilon G})\le\frac{\bar\nu_\phi}{8N_\phi}\,\varepsilon^2 .
\]
If $\phi$ is Lipschitz with constant $\mathrm{Lip}(\phi)$, then also $\omega_\phi(C_{\varepsilon G})\le\mathrm{Lip}(\phi)\Var(G)\,\varepsilon^2/(60N_\phi)$. Consequently no measure $\omega_\phi$ coincides with a measure $\kappa_\psi$. In particular $\xi$ is not of the form $\kappa_\psi$, and $\Rrho$ is not of the form $\omega_\phi$.
\end{proposition}

\begin{proof}
Spectral risk measures are positively homogeneous and translation equivariant, so $\mathcal M^+_\psi(L_t)-t=\varepsilon t(1-t)\mathcal M^+_\psi(G)$, and $\int_0^1t(1-t)\,\dd t=\tfrac16$; moreover $\mathcal M^+_\psi(G)>\E G=0$ because $\psi$ is not constant and $G$ is not degenerate. For $\omega_\phi$ put $Y:=L_t-t=\varepsilon t(1-t)\gamma(U)$, so that $\E Y=0$ and $|Y|\le a:=\varepsilon/4$, and extend $\phi$ to $\mathbb R$ by $\phi(0+)$ on $(-\infty,0]$ and $\phi(1-)$ on $[1,\infty)$. Since $\phi$ is nonincreasing,
\[
\Fp(t)-\Fp(t+Y)+\phi(t)Y=\int_0^Y\bigl(\phi(t)-\phi(t+s)\bigr)\,\dd s\le|Y|\bigl(\phi(t-a)-\phi(t+a)\bigr).
\]
Taking expectations and using $\E Y=0$ gives $\Fp(t)-\E\Fp(L_t)\le a\bigl(\phi(t-a)-\phi(t+a)\bigr)$, and $\int_0^1(\phi(t-a)-\phi(t+a))\,\dd t=\int_{-a}^a\phi-\int_{1-a}^{1+a}\phi\le2a\bar\nu_\phi$. This yields the bound $2a^2\bar\nu_\phi/N_\phi=\bar\nu_\phi\varepsilon^2/(8N_\phi)$. If $\phi$ is Lipschitz, the integral above is at most $\tfrac12\mathrm{Lip}(\phi) Y^2$, hence $\Fp(t)-\E\Fp(L_t)\le\tfrac12\mathrm{Lip}(\phi)\varepsilon^2t^2(1-t)^2\Var(G)$, and $\int_0^1t^2(1-t)^2\,\dd t=\tfrac1{30}$. Since $\kappa_\psi(C_{\varepsilon G})$ is linear in $\varepsilon$ with positive slope, it differs from $\omega_\phi(C_{\varepsilon G})$ for small $\varepsilon>0$.
\end{proof}

For $\phi=\phi_1$ one finds $\xi(C_{\varepsilon G})=\tfrac15\Var(G)\,\varepsilon^2$ and $\Rrho(C_{\varepsilon G})=\tfrac12\E|G-G'|\,\varepsilon$: along a ray, a scalar product is linear and a squared norm quadratic (Proposition~\ref{prop:pairing}). This is the dependence-theoretic form of the fact that an expected-utility functional is locally risk neutral, whereas a rank-dependent functional exhibits first-order risk aversion, and it explains why $\xi$ is typically much smaller than $\Rrho$ for weak dependence.

\begin{example}[$\xi$ and $\Rrho$ can rank copulas in opposite ways]\label{ex:lab}
In \eqref{eq:lab} with $\varepsilon=1$, let $G=\pm\tfrac12$ with probability $\tfrac12$ each. Then $(\xi,\Rrho)=(0.05,0.25)$. If instead $G\in\{-1,0,1\}$ with probabilities $(0.13,0.74,0.13)$, then $(\xi,\Rrho)=(0.052,0.226)$. Hence neither measure is a monotone function of the other, even on stochastically increasing copulas, where $\Rrho$ is Spearman's $\rho$, and the two copulas are not comparable in the Schur order. In fact $\Rrho$ is not a function of $\xi$ at all: for the second $G$ and $\varepsilon=(0.05/0.052)^{1/2}\approx0.981$ one obtains $\xi=0.05$, as for the first copula, but $\Rrho\approx0.222$.
\end{example}

\section{Sharp inequalities}\label{sec:sharp}

We now ask what the value of a spectral dependence measure implies for the values of the rearranged copula, that is, for the conditional probabilities $\Cr(x,t)/x$. By \eqref{eq:supform}, $\Cr(x,t)/x$ is the largest value of $\P(V\le t\mid U\in A)$ over all sets $A$ with $\P(U\in A)=x$; the supremum in \eqref{eq:supform} is attained by an indicator $\eta=\mathbf 1\{U\in A\}$, because $L_t$ is a function of $U$ and $U$ has a continuous distribution.

\subsection{Tent copulas and the main inequality}

For $x,t\in(0,1)$ and $0\le c\le M(x,t)-xt$ define the \emph{tent copula}
\begin{equation}\label{eq:tent}
T_{x,t,c}(u,v):=uv+A_x(u)\,g_{x,t,c}(v),
\end{equation}
where $A_x(u)=\min\{u/x,(1-u)/(1-x)\}$ and $g_{x,t,c}(v)=\bigl(c-(1-x)(t-v)_+-x(v-t)_+\bigr)_+$ is a tent function with peak $c$ at $v=t$, slope $1-x$ on the left and $-x$ on the right. We have $\partial_1T_{x,t,c}(u,v)=v+g_{x,t,c}(v)/x$ for $u<x$ and $v-g_{x,t,c}(v)/(1-x)$ for $u>x$. The derivatives of these functions with respect to $v$ take the values $1/x$, $0$, $1$ and $0$, $1/(1-x)$, $1$, respectively, so $T_{x,t,c}$ is a copula, and the condition on $c$ is exactly what is needed for the support of the tent to lie in $[0,1]$. Moreover $T_{x,t,c}$ is stochastically increasing, so it coincides with its rearrangement, and
\begin{equation}\label{eq:tenth}
(T_{x,t,c}-\Pi)(\alpha,v)=A_x(\alpha)\,g_{x,t,c}(v),\qquad\text{in particular}\quad(T_{x,t,c}-\Pi)(x,t)=c .
\end{equation}

For $b\in(0,1]$ we also consider the \emph{localised} measures
\[
\kappa_{\phi,b}(C):=\frac{1}{N_{\phi,b}}\int_0^b\bigl(\mathcal M^+_\phi(L_t)-t\bigr)\,\dd t,\qquad N_{\phi,b}:=\int_0^b\bigl(\Fp(t)-t\bigr)\,\dd t,
\]
so that $\kappa_{\phi,1}=\kappa_\phi$. They take values in $[0,1]$, are monotone in the Schur order, and vanish if and only if $L_t=t$ almost surely for all $t\le b$, that is, if and only if the events $\{V\le t\}$, $t\le b$, are independent of $U$.

\begin{theorem}\label{thm:sharp}
Let $C$ be a copula, $x,t\in(0,1)$, and $c:=(\Cr-\Pi)(x,t)$.
\begin{enumerate}
\item[\textup{(a)}] $\displaystyle\xi(C)\ge\frac{2\,c^3}{x^2(1-x)^2}$.
\item[\textup{(b)}] For every spectrum $\phi$ and every $b\in[t,1]$, with $\delta:=b-t$,
\[
N_{\phi,b}\,\kappa_{\phi,b}(C)\;\ge\;\frac{\Fp(x)-x}{x(1-x)}\Bigl(\frac{c^2}{2x(1-x)}-\frac x2\Bigl(\Bigl(\frac cx-\delta\Bigr)_+\Bigr)^2\Bigr).
\]
For $b=1$ the second term vanishes, so that
\[
\kappa_\phi(C)\ge\frac{(\Fp(x)-x)\,c^2}{2N_\phi\,x^2(1-x)^2},\qquad\text{in particular}\quad\Rrho(C)\ge\frac{3c^2}{x(1-x)} .
\]
\end{enumerate}
All inequalities are sharp for every admissible value of $c$: equality holds simultaneously in \textup{(a)} and, for all $\phi$ and $b$, in \textup{(b)} for the tent copula $T_{x,t,c}$.
\end{theorem}

\begin{proof}
\emph{Step 1 (in the direction of $x$).} Fix $v\in(0,1)$. Since $\partial_\alpha h_v(\alpha)=L_v^*(\alpha)-v$ has mean zero, $h_v(x)=\int_0^1(\mathbf 1_{[0,x]}(\alpha)-x)\,\partial_\alpha h_v(\alpha)\,\dd\alpha$, and the Cauchy--Schwarz inequality gives
\begin{equation}\label{eq:CS}
h_v(x)^2\le x(1-x)\int_0^1\bigl(\partial_\alpha h_v\bigr)^2\,\dd\alpha=x(1-x)\Var(L_v).
\end{equation}
Since $h_v$ is concave and vanishes at $0$ and $1$, its graph lies above the triangle with vertices $(0,0)$, $(x,h_v(x))$, $(1,0)$, that is, $h_v\ge h_v(x)A_x$. Hence, with $a_\phi(x):=\int A_x\,\dd\nu_\phi=\int_0^1\phi(s)A_x'(s)\,\dd s=(\Fp(x)-x)/(x(1-x))$,
\begin{equation}\label{eq:triangle}
\mathcal M^+_\phi(L_v)-v=\int_{(0,1)}h_v\,\dd\nu_\phi\ge a_\phi(x)\,h_v(x).
\end{equation}
\emph{Step 2 (in the direction of $t$).} By Lemma~\ref{lem:basic}(b) and (c), $v\mapsto h_v(x)$ is nonnegative with slopes in $[-x,1-x]$ and equals $c$ at $v=t$. Hence $h_v(x)\ge g_{x,t,c}(v)$ for all $v$.

\emph{Step 3 (integration).} We have $\int_0^1g_{x,t,c}^2=c^3/(3x(1-x))$, so Theorem~\ref{thm:dictionary}(i) and \eqref{eq:CS} give (a). The left part of the tent has base $c/(1-x)\le t$ and area $c^2/(2(1-x))$; the right part has base $c/x$, and its area over $[t,b]$ is $c^2/(2x)$ if $c/x\le\delta$ and $c\delta-x\delta^2/2$ otherwise. Thus $\int_0^bg_{x,t,c}$ equals the bracket in (b), and \eqref{eq:triangle} gives (b). For $b=1$ note that $c\le x(1-t)$. For $\phi=\phi_1$ we have $\Fp(x)-x=x(1-x)$ and $N_\phi=\tfrac16$.

\emph{Sharpness.} By \eqref{eq:tenth}, for $T_{x,t,c}$ the function $\partial_\alpha h_v$ is a multiple of $\mathbf 1_{[0,x]}-x$, the function $h_v$ is a multiple of $A_x$, and $h_v(x)=g_{x,t,c}(v)$. Hence \eqref{eq:CS}, \eqref{eq:triangle} and Step~2 are equalities.
\end{proof}

Inequality \eqref{eq:CS} is the sharp pointwise Sobolev inequality $|w(x)|^2\le x(1-x)\norm{w'}_{L^2}^2$ for functions $w$ on $[0,1]$ vanishing at the endpoints, whose extremal functions are the multiples of the Green's function $A_x$ of $-\dd^2/\dd x^2$; the tent copulas are the copulas for which $\Cr-\Pi$ is, for every $v$, a multiple of this representing function.

For $x=t=\tfrac12$ we have $c=\beta(\Cr)/4$ with Blomqvist's $\beta(C)=4C(\tfrac12,\tfrac12)-1$, and $|\beta(C)|\le\beta(\Cr)$ by Lemma~\ref{lem:basic}(f). Part (a) therefore becomes $\xi(C)\ge\beta(\Cr)^3/2\ge|\beta(C)|^3/2$. This is the inequality of \citet[Theorem~1.1]{OrendayRockel2026}, which describes the boundary of the exact $(\xi,\beta)$ region for all copulas and whose extremal copulas are the tent copulas $T_{1/2,1/2,c}$; our proof follows their argument, a Jensen-type estimate in the first variable followed by a one-dimensional minimisation over profiles with bounded slopes. We do not address the reverse inequalities, which bound $c$ from below in terms of the measures. For stochastically increasing copulas and $x=t=\tfrac12$ such an inequality is contained in \citet[Theorem~1.2]{OrendayRockel2026}.

\subsection{Consequences}

\begin{corollary}[Uniform distance from independence]\label{cor:uniform}
For every copula $C$,
\[
\norm{C-\Pi}_\infty\le\norm{\Cr-\Pi}_\infty\le\min\Bigl\{\Bigl(\frac{\xi(C)}{32}\Bigr)^{1/3},\Bigl(\frac{\Rrho(C)}{12}\Bigr)^{1/2}\Bigr\}.
\]
Both bounds are attained, for every value $\xi(C)\le\tfrac12$ and $\Rrho(C)\le\tfrac34$, respectively, by the copulas $T_{1/2,1/2,c}$, $0\le c\le\tfrac14$. Conversely,
\[
\xi(C)\le\zeta_1(C)\le G\bigl(\norm{\Cr-\Pi}_\infty\bigr),\qquad G(c):=6\int_0^1\min\{c,t(1-t)\}\,\dd t=6cs_c+6t_c^2-4t_c^3,
\]
where $s_c=\sqrt{1-4c}$ and $t_c=\tfrac12(1-s_c)$, and both inequalities are equalities for the copulas $Q_c$, $0\le c\le\tfrac14$, defined in the proof. In particular $G^{-1}(\xi(C))\le\norm{\Cr-\Pi}_\infty\le(\xi(C)/32)^{1/3}$, and both bounds are sharp.
\end{corollary}

\begin{proof}
The first inequality is Lemma~\ref{lem:basic}(f). Theorem~\ref{thm:sharp} gives $c\le(\xi/2)^{1/3}(x(1-x))^{2/3}$ and $c\le(\Rrho\,x(1-x)/3)^{1/2}$, and $x(1-x)\le\tfrac14$. For $T_{1/2,1/2,c}$ we have $\xi=32c^3$, $\Rrho=12c^2$ and $\norm{T-\Pi}_\infty=c$. For the converse let $c:=\norm{\Cr-\Pi}_\infty$ and $m_t:=\max_xh_t(x)$. Then $m_t\le c$ and, by Lemma~\ref{lem:basic}(b), $m_t\le t(1-t)$, so $\zeta_1=6\int_0^1m_t\,\dd t\le G(c)$ by Theorem~\ref{thm:dictionary}(iii), and $\xi\le\zeta_1$ by Corollary~\ref{cor:ineq}. Let $g_-(t):=\min\{c/(1-t),t\}$ and $g_+(t):=\min\{c/t,1-t\}$, and let $Q_c$ be given by $\partial_1Q_c(u,t)=1$ for $u\le g_-(t)$, $\partial_1Q_c(u,t)=t$ for $g_-(t)<u<1-g_+(t)$, and $\partial_1Q_c(u,t)=0$ for $u\ge1-g_+(t)$. Since $g_-$ is nondecreasing, $g_+$ is nonincreasing and $g_-(t)(1-t)=t\,g_+(t)=\min\{c,t(1-t)\}$, so that $g_-+g_+\le1$, this defines a stochastically increasing copula. Under $Q_c$ we have $L_t\in\{0,t,1\}$ and $m_t=\min\{c,t(1-t)\}$, hence $\norm{Q_c-\Pi}_\infty=c$ and $\xi(Q_c)=\zeta_1(Q_c)=G(c)$ by Corollary~\ref{cor:ineq}.
\end{proof}

The corollary is a sharp quantitative version of the statement that $\xi(C)=0$ implies independence, in terms of the supremum distance underlying the measure $\sigma_\infty$ of \citet{SchweizerWolff1981}. The lower bound concerns $\Cr$; no such bound holds for $\norm{C-\Pi}_\infty$, since the completely dependent copulas, for which $\xi=1$, are dense in the set of all copulas with respect to uniform convergence. Pointwise in $t$, the proofs of Corollary~\ref{cor:ineq} and of this corollary give $4\theta_t(t)^2\le\Var(L_t)\le\theta_t(t)$.

The situation is different for \emph{relative} deviations near the boundary. By Lemma~\ref{lem:basic}(b), the largest possible value of $\Cr(x,t)$ is $M(x,t)$.

\begin{corollary}[Global measures do not see the tail]\label{cor:tail}
Let $0<x\le t<1$ and $T:=T_{x,t,x(1-t)}$. Then $T(x,t)=M(x,t)=x$, that is, $\P(V\le t\mid U\le x)=1$, while
\[
\xi(T)=\frac{2x(1-t)^3}{(1-x)^2}<2x,\qquad\zeta_1(T)=\Rrho(T)=\frac{3x(1-t)^2}{1-x},
\]
and $\max\{\omega_\phi(T),\kappa_\phi(T)\}\le\bar\nu_\phi\,x/(2N_\phi(1-x))$ for every spectrum $\phi$. Conversely, $\Cr(x,t)=M(x,t)$ implies $\xi(C)\ge\xi(T)$ and $\kappa_\phi(C)\ge\kappa_\phi(T)$ for all spectra $\phi$.
\end{corollary}

\begin{proof}
Insert $c=x(1-t)$ into Theorem~\ref{thm:sharp} and use $\zeta_1(T)=6\int_0^1g_{x,t,c}=3c^2/(x(1-x))$, which follows from Theorem~\ref{thm:dictionary}(iii) and \eqref{eq:tenth}, together with Proposition~\ref{prop:cont}.
\end{proof}

In terms of the conditional probability $\Cr(x,t)/x$, Theorem~\ref{thm:sharp} gives the bounds
\begin{equation}\label{eq:blind}
\frac{\Cr(x,t)}{x}-t\;\le\;\min\Bigl\{\Bigl(\frac{\xi(C)}{2}\Bigr)^{1/3}\frac{(1-x)^{2/3}}{x^{1/3}},\ \Bigl(\frac{\Rrho(C)\,(1-x)}{3x}\Bigr)^{1/2}\Bigr\},
\end{equation}
which are sharp as long as they do not exceed the trivial bound $M(x,t)/x-t$, and which diverge as $x\to0$. To control the tail, the spectrum has to be adapted to the tail level.

\begin{corollary}[Tail spectra]\label{cor:stress}
For $p\in(0,1]$ let $\phi_p(s)=2(p-s)_+/p^2$, so that $\phi_1$ is the spectrum of Section~\ref{sec:spectral}, and, for $b\in[p,1]$, $\tau_{p,b}:=\kappa_{\phi_p,b}$ and $\tau_p:=\tau_{p,1}$. Then $N_{\phi_p,b}=\tfrac16p(4-3p)+\tfrac12(b-p)(2-b-p)$, and if $p<1$, then for every copula $C$ and $t\in(0,b]$, with $e:=\Cr(p,t)/p-t$ and $\delta:=b-t$,
\begin{equation}\label{eq:stress}
N_{\phi_p,b}\,\tau_{p,b}(C)\;\ge\;\frac{e^2}{2(1-p)}-\frac{\bigl((e-\delta)_+\bigr)^2}{2},\qquad\text{in particular}\quad e\le\sqrt{\tfrac13(1-p)(3-2p)\,\tau_p(C)}.
\end{equation}
Equality holds for $T_{p,t,pe}$. Moreover $\tau_{1,1}=\Rrho$, and $\tau_{p,b}(C)=1$ if and only if $L_t\in\{0,1\}$ almost surely for all $t\le p$.
\end{corollary}

\begin{proof}
Here $F_{\phi_p}(s)=1-(p-s)_+^2/p^2$, which gives $N_{\phi_p,b}$ and $F_{\phi_p}(p)-p=1-p$. Inequality \eqref{eq:stress} is Theorem~\ref{thm:sharp}(b) with $x=p$ and $c=pe$; for $b=1$ we have $e\le1-t=\delta$ and $N_{\phi_p,1}=(3-2p)/6$. If $\tau_{p,b}(C)=1$, then $h_t=h^M_t$ on $[0,p]$ for almost every $t\le b$, because the support of $\nu_{\phi_p}$ is $[0,p]$; for almost every $t\le p$ the choice $x=t$ gives $L_t\in\{0,1\}$, and this extends to all $t\le p$ because $t\mapsto L_t$ is continuous in $L^1$ and the set of $\{0,1\}$-valued random variables is closed in $L^1$. Conversely, if $L_t\in\{0,1\}$ for all $t\le p$, then for $t\in(p,b]$ we have $L_t\ge L_p=1$ on a set of probability $p$, so $h_t(\alpha)=\alpha(1-t)=h^M_t(\alpha)$ for $\alpha\le p$.
\end{proof}

The right-hand side of the first inequality in \eqref{eq:stress} is increasing in $e$, so the inequality is an upper bound for $e$ in terms of $\tau_{p,b}(C)$. In contrast to \eqref{eq:blind}, it does not degenerate as $p\to0$: for $b=1$ it reads, in the limit, $\operatorname{ess\,sup}L_t-t\le\sqrt{\tau_0(C)}$, where $\tau_0(C)=2\int_0^1(\operatorname{ess\,sup}L_t-t)\,\dd t$. Like all measures $\kappa_{\phi,b}$ they are convex in $C$, depend on $C$ only through $\Cr$, and satisfy the information inequality of Theorem~\ref{thm:mirror}(a), with the same proofs. For the Farlie--Gumbel--Morgenstern copula with parameter $\alpha\in[0,1]$ one finds $\tau_p=\alpha/3$ for all $p$, so in this family the tail level is irrelevant; in general it is not, as Section~\ref{sec:illustration} shows.

\subsection{The stop-loss diagram}\label{sec:stoploss}

Theorem~\ref{thm:sharp} concerns the Lorenz diagram. We now give the corresponding statements for the stop-loss excess $\theta_t(q)$, whose largest possible value is $M(q,t)-qt$ by Lemma~\ref{lem:basic}(d), and for the indicator correlation $\varrho_C(t)$, which localises $\xi$ in $t$. Let $B_t:=T_{t,t,t(1-t)}$ denote the \emph{block copula}, under which $V$ is uniformly distributed on $[0,t]$ if $U<t$ and on $[t,1]$ otherwise.

\begin{proposition}[Pointwise bounds]\label{prop:profile}
For every copula $C$ and all $x,q,t\in(0,1)$,
\begin{align}
(\Cr-\Pi)(x,t)&\le\sqrt{x(1-x)\,t(1-t)\,\varrho_C(t)},\label{eq:profES}\\
\theta_t(q)&\le\tfrac12\Bigl(\sqrt{t(1-t)\,\varrho_C(t)+(q-t)^2}-|q-t|\Bigr).\label{eq:profTranche}
\end{align}
Equality holds in \eqref{eq:profES} if and only if $L_t$ takes at most two values, the larger one with probability $x$. Equality holds in \eqref{eq:profTranche} if and only if $L_t$ takes values in $\{q-r,q+r\}$ for some $r\ge0$; such a law with mean $t$ and variance $t(1-t)\varrho_C(t)$ exists only if its support lies in $[0,1]$. Moreover, if the limit exists, $\lim_{t\to0}\varrho_C(t)$ is the lower tail dependence coefficient of the copula of $(V,V')$.
\end{proposition}

\begin{proof}
\eqref{eq:profES} is \eqref{eq:CS}, in which equality holds if and only if $L_t^*-t$ is a multiple of $\mathbf 1_{[0,x]}-x$. For a random variable $X$ with mean $\mu$ and variance $\sigma^2$ we have $(X-q)_+=\tfrac12(|X-q|+X-q)$ and $\E|X-q|\le\sqrt{\sigma^2+(\mu-q)^2}$, hence $\E(X-q)_+-(\mu-q)_+\le\tfrac12\bigl(\sqrt{\sigma^2+(q-\mu)^2}-|q-\mu|\bigr)$; this is the bound of \citet{Scarf1958} and \citet{Lo1987}, with equality if and only if $|X-q|$ is constant. Finally, $\varrho_C(t)=(\E L_t^2/t-t)/(1-t)$ and $\E L_t^2/t=\P(V'\le t\mid V\le t)$.
\end{proof}

\begin{theorem}\label{thm:stoploss}
Let $C$ be a copula, $q,t\in(0,1)$, $\theta:=\theta_t(q)$, and $\delta(x):=(t-q)_+-x(t-q)$ for $x\in[0,1]$.
\begin{enumerate}
\item[\textup{(a)}] If $\theta>0$, then
\begin{equation}\label{eq:stoplossSharp}
\xi(C)\;\ge\;\Xi_{t,q}(\theta):=\min\Bigl\{\frac{2\,(\theta+\delta(x))^3}{x^2(1-x)^2}:\ x\in(0,1),\ \theta+\delta(x)\le M(x,t)-xt\Bigr\},
\end{equation}
where the admissible $x$ form a nonempty compact interval. Equality holds for the tent copula $T_{\bar x,t,\theta+\delta(\bar x)}$, where $\bar x$ is a minimiser; it is the point of the admissible interval closest to $x_0$ if $q\ge t$ and to $1-x_0$ if $q<t$, where $x_0=\bigl(\sqrt{(d+4\theta)^2+8\theta d}-d-4\theta\bigr)/(2d)$ with $d:=|q-t|$, and $x_0=\tfrac12$ if $d=0$.
\item[\textup{(b)}] $\theta_t(q)=M(q,t)-qt$ holds if and only if $L_t\in\{0,1\}$ almost surely, and then it holds for every $q\in(0,1)$. For all such copulas $\xi(C)\ge2t(1-t)$, with equality if and only if $\Cr=B_t$.
\end{enumerate}
\end{theorem}

\begin{proof}
(a) By Lemma~\ref{lem:basic}(e), $\theta=\max_x\{S_t(x)-xq\}-(t-q)_+=\max_x\{h_t(x)-\delta(x)\}$. The function $h_t-\delta$ is concave, with values $-(t-q)_+\le0$ at $0$ and $-(q-t)_+\le0$ at $1$, so its maximum $\theta>0$ is attained at some $x^*\in(0,1)$, and $c:=h_t(x^*)=\theta+\delta(x^*)$. By Lemma~\ref{lem:basic}(b), $c\le M(x^*,t)-x^*t$, so $x^*$ is admissible, and Theorem~\ref{thm:sharp}(a) gives $\xi(C)\ge2c^3/(x^{*2}(1-x^*)^2)\ge\Xi_{t,q}(\theta)$. Both admissibility constraints are linear in $x$; for $q\ge t$ they describe the interval $[\theta/(1-q),(t-\theta)/q]$ and for $q<t$ the interval $[(\theta+t-q)/(1-q),1-\theta/q]$, which contain $x^*$ and lie in $(0,1)$. Conversely, for admissible $x$ let $T:=T_{x,t,\theta+\delta(x)}$. By \eqref{eq:tenth}, $h_t^T-\delta=(\theta+\delta(x))A_x-\delta$ is concave and piecewise linear with a single kink at $x$, where it takes the value $\theta$, and it is nonpositive at $0$ and $1$; hence $\theta_t^T(q)=\theta$, while $\xi(T)=2(\theta+\delta(x))^3/(x^2(1-x)^2)$ by Theorem~\ref{thm:sharp}. This proves \eqref{eq:stoplossSharp} and the equality statement. Finally, for $q\ge t$ the logarithmic derivative of $2(\theta+xd)^3/(x^2(1-x)^2)$ vanishes if and only if $dx^2+(d+4\theta)x-2\theta=0$, which has the single root $x_0$ in $(0,1)$; the function decreases to the left and increases to the right of $x_0$, so its minimum over an interval is attained at the point closest to $x_0$. For $q<t$ replace $x$ by $1-x$.

(b) For $q\ge t$ the equality means $\E(L_t-q)_+=(1-q)\,\E L_t$, and $(\ell-q)_+<(1-q)\ell$ for $\ell\in(0,1)$; for $q<t$ one argues with $\E(q-L_t)_+=q\,\E(1-L_t)$. Hence $L_t=\mathbf 1_A$ with $\P(A)=t$. For $v<t$ monotonicity gives $L_v=0$ outside $A$, so $\E[L_v\mid A]=v/t$ and $\Var(L_v)\ge v^2/t-v^2$, with equality if and only if $L_v=v/t$ on $A$. For $v>t$ we have $L_v=1$ on $A$ and $\E[L_v\mid A^c]=(v-t)/(1-t)$, so $\Var(L_v)\ge t(1-v)^2/(1-t)$, with equality if and only if $L_v$ is constant on $A^c$. Integration gives $\int_0^1\Var(L_v)\,\dd v\ge\tfrac13t(1-t)$, with equality if and only if $L_v$ has the same law as under $B_t$ for almost every $v$, hence, by continuity of $v\mapsto L_v$ in $L^1$, for every $v$; that is, $\Cr=B_t^\uparrow=B_t$.
\end{proof}

Part (b) is the stop-loss counterpart of Corollary~\ref{cor:tail} with $x=t$. Part (a) determines the smallest value of $\xi$ that is compatible with a given stop-loss excess $\theta_t(q)$; it therefore dominates every other lower bound of this kind, for instance the bound $\xi\ge12\theta(\theta+d)(2\theta+d)$ that follows from the Scarf--Lo bound in Proposition~\ref{prop:profile} and Lemma~\ref{lem:basic}(g). For $q=t$ and $\theta\le\tfrac12\min\{t,1-t\}$ it reads $\xi\ge32\theta^3$, in accordance with Corollary~\ref{cor:uniform}, since $\theta_t(t)=\max_xh_t(x)$. By Theorem~\ref{thm:dictionary}(i), $\xi$ is the average of the profile $\varrho_C$ with weight $6t(1-t)$, which vanishes at the boundary. This is a second reason why $\xi$ does not see the lower tail: for the Gaussian copula with correlation $r$ the pair $(V,V')$ has a Gaussian copula with correlation $r^2$, which has no tail dependence, so $\varrho_C(t)\to0$ as $t\to0$ whatever the value of $\xi$.

\section{Illustration}\label{sec:illustration}

We compare the Gaussian copula with correlation $r=\sqrt{0.2}$, for which $\xi=\tfrac3\pi\arcsin\bigl(\tfrac{1+r^2}{2}\bigr)-\tfrac12=0.1145$, with the Student $t$ copula with four degrees of freedom and correlation parameter $0.441$, chosen such that $\xi$ takes the same value. The Gaussian copula is stochastically increasing, so that $\Cr=C$ and $\Rrho=\rho=\tfrac6\pi\arcsin(r/2)$; the $t$ copula is not, and $\Cr\ne C$. We add the tent copula of Corollary~\ref{cor:tail} with $x=0.001$ and $t=0.01$ and the block copula $B_t$ of Section~\ref{sec:stoploss} with $t=0.01$, for which all quantities are available in closed form. For the two parametric copulas the values of $\partial_1C(\cdot,t)$ were sorted on a grid of $6\cdot10^4$ points that is uniform on the normal quantile scale, and spectral functionals were evaluated by exact integration of the spectrum over the cells; integrals over $t$ were computed by composite Gauss--Legendre rules. This reproduces the closed-form values for the Farlie--Gumbel--Morgenstern and the Gaussian copula to at least six digits, and for the Student copula it agrees to six digits with an evaluation of the variational formula of Lemma~\ref{lem:basic}(e) by adaptive quadrature.

\begin{table}[htbp]
\centering
\caption{Global and tail-localised dependence measures; $p=0.001$.}\label{tab:global}
\begin{tabular}{lccccc}
\toprule
Copula & $\xi$ & $\Rrho=\tau_{1,1}$ & $\zeta_1$ & $\tau_{p,1}$ & $\tau_{p,0.05}$\\
\midrule
Gaussian & 0.1145 & 0.431 & 0.312 & 0.755 & 0.292\\
Student $t_4$ & 0.1145 & 0.414 & 0.312 & 0.679 & 0.610\\
Tent $T_{x,t,x(1-t)}$, $x=0.001$, $t=0.01$ & 0.0019 & 0.003 & 0.003 & 0.982 & 0.812\\
Block $B_t$, $t=0.01$ & 0.0198 & 0.030 & 0.030 & 0.991 & 0.904\\
\bottomrule
\end{tabular}
\end{table}

\begin{table}[htbp]
\centering
\caption{The conditional probability $\Cr(x,t)/x$ at $x=0.001$, and upper bounds on $\Cr(x,t)/x-t$ at $t=0.01$ implied by single quantities via \eqref{eq:blind}, \eqref{eq:stress} and \eqref{eq:profES}, with $p=x$. The trivial bound is $1-t=0.99$.}\label{tab:local}
\begin{tabular}{lcccccccc}
\toprule
 & \multicolumn{3}{c}{$\Cr(x,t)/x$} & \multicolumn{5}{c}{bound on $\Cr(x,t)/x-t$ at $t=0.01$ from}\\
\cmidrule(lr){2-4}\cmidrule(lr){5-9}
 & $t=0.001$ & $t=0.01$ & $t=0.05$ & $\xi$ & $\Rrho$ & $\tau_{p,1}$ & $\tau_{p,0.05}$ & $\varrho_C(t)$\\
\midrule
Gaussian & 0.040 & 0.181 & 0.439 & 0.99 & 0.99 & 0.868 & 0.372 & 0.489\\
Student $t_4$ & 0.232 & 0.524 & 0.680 & 0.99 & 0.99 & 0.823 & 0.751 & 0.959\\
\bottomrule
\end{tabular}
\end{table}

\begin{figure}[htbp]
\centering
\includegraphics[width=\textwidth]{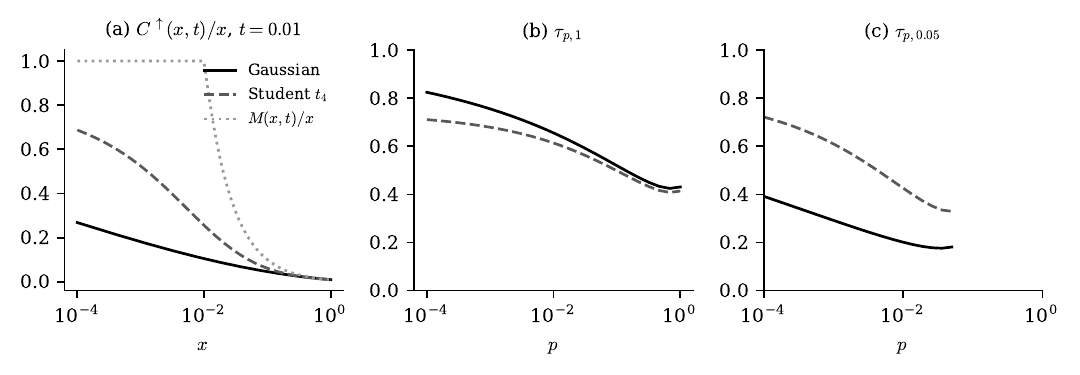}
\caption{Gaussian and Student copula with identical $\xi=0.1145$. (a) The conditional probability $x\mapsto\Cr(x,t)/x$ at $t=0.01$, together with its largest possible value $M(x,t)/x$. (b), (c) The tail measures $\tau_{p,b}$ as functions of $p\le b$ for $b=1$ and $b=0.05$.}\label{fig:tail}
\end{figure}

Table~\ref{tab:global} shows that the two parametric copulas are practically indistinguishable for the global measures: $\xi$ is equal by construction, $\zeta_1$ agrees to three digits, and $\Rrho$ is slightly smaller for the Student copula. The rearranged copulas are nevertheless very different near the origin (Table~\ref{tab:local} and Figure~\ref{fig:tail}(a)): the largest conditional probability of $\{V\le0.01\}$ given an event of probability $0.001$ in $U$ is $0.18$ for the Gaussian and $0.52$ for the Student copula. In accordance with \eqref{eq:blind}, the bounds implied by $\xi$ and $\Rrho$ are trivial at this point. The tail measure $\tau_{p,1}$ gives a nontrivial but loose bound and still ranks the Gaussian above the Student copula, because it averages over all $t$ and the Gaussian copula has the larger excess for $t$ near $\tfrac12$. Localising to $t\le0.05$ reverses the ranking and separates the two copulas: $\tau_{0.001,0.05}$ equals $0.29$ and $0.61$, and the bounds from \eqref{eq:stress} are $0.37$ and $0.75$, against true values of $0.17$ and $0.51$. Since \eqref{eq:stress} is sharp, the remaining gap is attained by a tent copula with the same value of $\tau_{0.001,0.05}$. The two extremal copulas in the last rows of Table~\ref{tab:global} have negligible global measures and tail measures close to one.

\section{Estimation}\label{sec:estimation}

All measures of this paper are functionals of the copula and can be estimated by evaluating them at an estimator of $C$ that is itself a copula. By Proposition~\ref{prop:cont}, Trutschnig's metric $D_1$ is the natural metric for this purpose. We use the empirical checkerboard copula, for which $D_1$-consistency and, under a smoothness assumption, a rate of convergence are known.

\subsection{The estimator}

Let $(X_1,Y_1),\dots,(X_N,Y_N)$ be independent copies of a random vector with continuous marginal distribution functions and copula $C$, and let $R_k$ and $S_k$ denote the ranks of $X_k$ among $X_1,\dots,X_N$ and of $Y_k$ among $Y_1,\dots,Y_N$. Following \citet[Section~4.3]{AnsariWiesel2026}, we assume $N=m^3$ with an integer resolution $m$ and put
\[
a_{ij}:=\frac1{m^2}\,\#\bigl\{k:\ (i-1)m^2<R_k\le im^2,\ (j-1)m^2<S_k\le jm^2\bigr\},\qquad i,j=1,\dots,m .
\]
Every row and every column of $A_N=(a_{ij})$ sums to one, and the \emph{empirical checkerboard copula} $\hat C_N$ is the copula with density $m\,a_{ij}$ on the cell $\bigl(\frac{i-1}m,\frac im\bigr]\times\bigl(\frac{j-1}m,\frac jm\bigr]$. Thus $\hat C_N$ is the $m$-checkerboard approximation of the empirical copula, that is, the empirical checkerboard estimator of \citet{SDS2024} with resolution $N^{1/3}$ (they allow resolutions $\lfloor N^{s}\rfloor$ with $s\in(0,\tfrac12)$); the restriction $N=m^3$ only serves to make all row and column sums exactly one. Since $\hat C_N$ is a copula, $\omega_\phi(\hat C_N)$, $\kappa_{\phi,b}(\hat C_N)$ and in particular $\tau_{p,b}(\hat C_N)$ are well defined.

\begin{remark}[Exact evaluation]\label{rem:compute}
Let $D$ be a checkerboard copula with resolution $m$ and doubly stochastic matrix $A=(a_{ij})$ as above. For $u$ in the $i$-th strip $\bigl(\frac{i-1}m,\frac im\bigr]$ and $t=\frac{k-1+s}m$ with $s\in[0,1]$,
\[
\partial_1D(u,t)=\alpha_{ik}+a_{ik}\,s,\qquad \alpha_{ik}:=\sum_{j<k}a_{ij}.
\]
Hence $L_t$ takes the values $\ell_i(t)=\alpha_{ik}+a_{ik}s$, $i=1,\dots,m$, each with probability $\frac1m$, and
\[
\mathcal M^+_\phi(L_t)=\sum_{r=1}^m\ell_{(r)}(t)\Bigl(\Fp\bigl(\tfrac rm\bigr)-\Fp\bigl(\tfrac{r-1}m\bigr)\Bigr),\qquad \E\Fp(L_t)=\frac1m\sum_{i=1}^m\Fp\bigl(\ell_i(t)\bigr),
\]
where $\ell_{(1)}(t)\ge\dots\ge\ell_{(m)}(t)$. In every strip of $t$ the order of the affine functions $\ell_i$ changes only at their at most $\binom m2$ intersection points. Between consecutive intersection points, points where some $\ell_i(t)=t$, and the point $t=b$, the integrands of $\xi$, $\zeta_1$ and $\kappa_{\phi,b}$ are polynomials of degree at most two in $t$, so these measures are computed exactly by the two-point Gauss rule on each subinterval, with $O(m^4\log m)=O(N^{4/3}\log N)$ operations. No rearrangement of $D$ has to be formed explicitly.
\end{remark}

\subsection{Consistency and rates}

The proof of Proposition~\ref{prop:cont} applies verbatim to the localised measures and gives, for every spectrum $\phi$, $b\in(0,1]$ and all copulas $C,D$,
\begin{equation}\label{eq:contloc}
\bigl|\kappa_{\phi,b}(C)-\kappa_{\phi,b}(D)\bigr|\le\frac{\bar\nu_\phi}{2N_{\phi,b}}\int_0^b\!\!\int_0^1\bigl|\partial_1C(u,t)-\partial_1D(u,t)\bigr|\,\dd u\,\dd t\le\frac{\bar\nu_\phi}{2N_{\phi,b}}\,D_1(C,D).
\end{equation}
We say that $C$ has \emph{Lipschitz kernels} if there is $\Lambda<\infty$ with
\[
W_1\bigl(K_C(u,\cdot),K_C(u',\cdot)\bigr)=\int_0^1\bigl|\partial_1C(u,v)-\partial_1C(u',v)\bigr|\,\dd v\le\Lambda\,|u-u'|\qquad\text{for almost all }u,u',
\]
where $W_1$ is the Wasserstein distance on $[0,1]$; this is the assumption of \citet[Definition~4.1 and Theorem~4.7]{AnsariWiesel2026}, applied to the copula. It holds, for instance, if $C$ has a density $c$ with $\sup_u\int_0^1|\partial_uc(u,v)|\,\dd v<\infty$, since $\partial_u\partial_1C(u,v)=\int_0^v\partial_uc(u,y)\,\dd y$; for the Farlie--Gumbel--Morgenstern copula with parameter $\alpha$ one may take $\Lambda=|\alpha|$.

\begin{theorem}\label{thm:estimation}
Let $N=m^3$ and let $\hat C_N$ be the empirical checkerboard copula. For $K>0$ let $\mathcal S_K$ be the set of pairs $(\phi,b)$ of a spectrum $\phi$ and $b\in(0,1]$ with $\bar\nu_\phi\le K\,N_{\phi,b}$.
\begin{itemize}
\item[\textup{(a)}] Almost surely, as $m\to\infty$, for every $K>0$,
\[
\sup_{(\phi,b)\in\mathcal S_K}\bigl|\kappa_{\phi,b}(\hat C_N)-\kappa_{\phi,b}(C)\bigr|\to0,\qquad \sup_{(\phi,1)\in\mathcal S_K}\bigl|\omega_\phi(\hat C_N)-\omega_\phi(C)\bigr|\to0 .
\]
In particular $\sup\{|\tau_{p,b}(\hat C_N)-\tau_{p,b}(C)|:\ p_0\le p\le b\le1\}\to0$ almost surely for every $p_0>0$.
\item[\textup{(b)}] If $C$ has Lipschitz kernels, there is a constant $c$ such that for all $m$ and $K$
\[
\E\Bigl[\sup_{(\phi,b)\in\mathcal S_K}\bigl|\kappa_{\phi,b}(\hat C_N)-\kappa_{\phi,b}(C)\bigr|\Bigr]+\E\Bigl[\sup_{(\phi,1)\in\mathcal S_K}\bigl|\omega_\phi(\hat C_N)-\omega_\phi(C)\bigr|\Bigr]\le K\,c\,N^{-1/3}.
\]
\item[\textup{(c)}] If $C=M$, then $\hat C_N$ is the checkerboard copula $D_m$ with $A_N$ the identity matrix, for every sample. For every spectrum $\phi$ and $b\in(0,1]$,
\[
\omega_\phi(D_m)=1-\frac1m,\qquad \zeta_1(D_m)=1-\frac1{2m},\qquad 0\le1-\kappa_{\phi,b}(D_m)\le\frac{\bar\nu_\phi}{6N_{\phi,b}\,m^2},
\]
and $\Rrho(D_m)=1-m^{-2}$. Since $M$ has Lipschitz kernels with $\Lambda=1$, the rate $N^{-1/3}$ in \textup{(b)} cannot be improved for the measures $\omega_\phi$, and neither can the rate of $\E D_1(\hat C_N,C)$ in \citet[Theorem~4.7]{AnsariWiesel2026}.
\end{itemize}
\end{theorem}

\begin{proof}
(a) and (b) follow from Proposition~\ref{prop:cont} and \eqref{eq:contloc}, because the right-hand sides are at most $\tfrac K2D_1(\hat C_N,C)$ uniformly on $\mathcal S_K$, together with $D_1(\hat C_N,C)\to0$ almost surely \citep[Lemma~4.5]{AnsariWiesel2026} and, under the assumption of (b), $\E D_1(\hat C_N,C)\le c'N^{-1/3}$ \citep[Theorem~4.7]{AnsariWiesel2026}; note that their metric $d_1$ is $D_1$ and that their kernel condition refers to the copula, because $(F_X,F_Y)_\#\P^{(X,Y)}$ has distribution function $C$. For the tail spectra, $\bar\nu_{\phi_p}=2/p$ and $N_{\phi_p,b}\ge N_{\phi_p,p}=\tfrac16p(4-3p)\ge\tfrac p6$ by Corollary~\ref{cor:stress}, so $(\phi_p,b)\in\mathcal S_K$ with $K=12/p_0^2$ whenever $p_0\le p\le b$.

(c) If $C=M$, then $R_k=S_k$ for all $k$, so $A_N$ is the identity. For $t=\frac{k-1+s}m$ the variable $L_t$ equals one on the first $k-1$ strips, $s$ on the $k$-th strip and zero otherwise. Hence $\E\Fp(L_t)=\frac1m(k-1+\Fp(s))$ and $\int_0^1\E\Fp(L_t)\,\dd t=\frac{m-1}{2m}+\frac1m\int_0^1\Fp$; with $N_\phi=\int_0^1\Fp-\tfrac12$ and \eqref{eq:xiphi} this gives $N_\phi\,\omega_\phi(D_m)=(1-\frac1m)N_\phi$. The value of $\zeta_1(D_m)$ follows by integrating $\E|L_t-t|=\frac1m\bigl((k-1)(1-t)+|s-t|+(m-k)t\bigr)$. For $\kappa_{\phi,b}$, the Lorenz curves of $L_t$ under $M$ and under $D_m$ differ by the tent function $g_t$ that vanishes outside $I_k=(\frac{k-1}m,\frac km)$, is affine on $[\frac{k-1}m,t]$ and on $[t,\frac km]$, and satisfies $g_t(t)=s(1-s)/m$. Since $\mathcal M^+_\phi(\mathbf 1_{\{U\le t\}})=\Fp(t)$, we have $\kappa_{\phi,b}(M)=1$, and \eqref{eq:kappaphi}, localised to $[0,b]$, gives
\[
N_{\phi,b}\bigl(1-\kappa_{\phi,b}(D_m)\bigr)=\int_0^b\!\!\int g_t\,\dd\nu_\phi\,\dd t\le\sum_{k=1}^m\nu_\phi(I_k)\int_{I_k}\frac{s(1-s)}m\,\dd t=\frac{\bar\nu_\phi}{6m^2}.
\]
For $\phi_1$, $\nu_{\phi_1}$ is twice Lebesgue measure, so $\int g_t\,\dd\nu_{\phi_1}=s(1-s)/m^2$ and $1-\Rrho(D_m)=6\cdot\frac1{6m^2}$. Finally, $W_1(K_M(u,\cdot),K_M(u',\cdot))=|u-u'|$, and $|\xi(D_m)-\xi(M)|=\frac1m=N^{-1/3}$ together with Proposition~\ref{prop:cont} gives $D_1(D_m,M)\ge\frac1{6m}$.
\end{proof}

By part (c), the smoothing by the checkerboard affects $\omega_\phi$ at the order $1/m$ of the strip width, but $\kappa_\phi$ only at the order $1/m^2$, because rearrangement and integration against $\nu_\phi$ average the error over each strip. We do not know whether $N^{-1/3}$ is the optimal rate for $\kappa_\phi$ under the assumption of (b); \citet{AnsariWiesel2026} conjecture rates between $N^{-1/3}$ and $N^{-2/3}$ for rearranged measures.

The assumption of (b) excludes the models of Section~\ref{sec:illustration}. For the Gaussian copula with correlation $r\in(0,1)$, $z=\Phi^{-1}(u)$ and $\sigma=\sqrt{1-r^2}$, we have $\partial_1C(u,v)=\Phi\bigl((\Phi^{-1}(v)-rz)/\sigma\bigr)$, so $\partial_u\partial_1C(u,v)\le0$ for all $v$, and the identity, in which $\varphi$ denotes the standard normal density, $\int\varphi\bigl(\frac{y-a}\sigma\bigr)\varphi(y)\,\dd y=\frac{\sigma}{\sqrt{1+\sigma^2}}\,\varphi\bigl(\frac a{\sqrt{1+\sigma^2}}\bigr)$ gives
\[
W_1\bigl(K_C(u,\cdot),K_C(u',\cdot)\bigr)=\Bigl|\int_u^{u'}g\Bigr|,\qquad g(w)=\frac r{\sqrt{2-r^2}}\exp\Bigl(\frac{1-r^2}{2-r^2}\,\Phi^{-1}(w)^2\Bigr),
\]
which is unbounded near $w=0$ and $w=1$; numerically the same holds for the Student copula. Theorem~\ref{thm:estimation}(a) applies to both models, but no rate is available for them. The simulations below indicate that the rate is in fact better than $N^{-1/3}$.

\subsection{The resolution limit}

Theorem~\ref{thm:estimation} is uniform only over spectra with $\bar\nu_\phi\le KN_{\phi,b}$, and for the tail spectra this excludes $p\to0$. The following observation shows that this is not an artefact of the proof.

\begin{proposition}\label{prop:resolution}
Let $D$ be a checkerboard copula with resolution $m$, and let $b\in(0,1]$. Then
\[
N_{\phi_p,b}\,\tau_{p,b}(D)=\int_0^b\Bigl(\max_{u}\partial_1D(u,t)-t\Bigr)\,\dd t\qquad\text{for all }p\in\bigl(0,\min\{b,\tfrac1m\}\bigr].
\]
\end{proposition}

\begin{proof}
By Remark~\ref{rem:compute} the decreasing rearrangement $L_t^*$ is constant on $[0,\frac1m)$ with value $\max_u\partial_1D(u,t)$, so $h_t(x)=x\,(\max_u\partial_1D(u,t)-t)$ for $x\le\frac1m$. Since $\nu_{\phi_p}$ is $2p^{-2}$ times Lebesgue measure on $(0,p)$, we get $\int h_t\,\dd\nu_{\phi_p}=\max_u\partial_1D(u,t)-t$ for $p\le\frac1m$, and \eqref{eq:kappaphi}, localised to $[0,b]$, gives the claim.
\end{proof}

Hence, for $p<N^{-1/3}$, the estimate $\tau_{p,b}(\hat C_N)$ depends on $p$ only through the normalising constant and carries no information about the tail level $p$: what is estimated is the average largest conditional probability over the strips of width $N^{-1/3}$. With $N=10^4$ this restricts the resolvable tail levels to $p\gtrsim0.05$. The measures $\tau_{0.001,0.05}$ of Section~\ref{sec:illustration} would require $m\ge1000$, that is, $N\ge10^9$ observations, and are therefore population quantities that cannot be estimated from samples of realistic size without further modelling assumptions on the tail.

\subsection{A simulation study}\label{sec:simulation}

\begin{figure}[htbp]
\centering
\includegraphics[width=\textwidth]{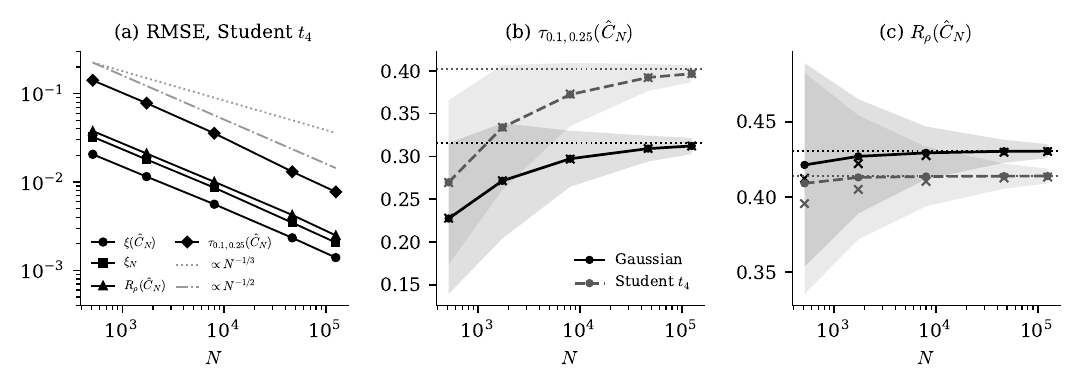}
\caption{(a) Root mean squared error of four estimators for the Student copula, with reference slopes. (b), (c) Monte Carlo means of $\tau_{0.1,0.25}(\hat C_N)$ and $\Rrho(\hat C_N)$ with pointwise $95\%$ ranges (mean $\pm1.96$ standard deviations, shaded), population values (dotted) and values for the checkerboard approximation $C_m$ (crosses).}\label{fig:est}
\end{figure}

We use the Gaussian and the Student copula of Section~\ref{sec:illustration}, which have the same value $\xi=0.1145$, and the sample sizes $N=m^3$ with $m\in\{8,12,20,36,50\}$, that is, $N$ between $512$ and $125\,000$. For each model and sample size we draw $1000$ samples and compute $\xi(\hat C_N)$, $\Rrho(\hat C_N)$, $\tau_{0.25,0.25}(\hat C_N)$ and $\tau_{0.1,0.25}(\hat C_N)$ exactly as in Remark~\ref{rem:compute}, and, for comparison, Chatterjee's rank correlation $\xi_N=1-3\sum_{k=1}^{N-1}|r_{k+1}-r_k|/(N^2-1)$, where $r_k$ is the rank of the $Y$-value of the observation with the $k$-th smallest $X$-value. The tail levels are chosen with Proposition~\ref{prop:resolution} in mind: $p=0.1$ exceeds $1/m$ except for $m=8$. The population values, computed as in Section~\ref{sec:illustration}, are given in Table~\ref{tab:sim}; the two localised measures rank the copulas oppositely to $\Rrho$. To separate the effect of the resolution from sampling error we also computed the measures for the $m$-checkerboard approximation $C_m$ of $C$ itself, whose cell masses are the $C$-volumes of the cells; we call $\kappa(C_m)-\kappa(C)$ the \emph{resolution bias}. The results are summarised in Table~\ref{tab:sim} and Figure~\ref{fig:est}.

\begin{table}[htbp]
\centering
\caption{Monte Carlo means and standard deviations (in parentheses) of the plug-in estimators, $1000$ replications per entry, together with the resolution bias for the Student copula and the probability that the estimate from a Student sample exceeds the estimate from an independent Gaussian sample of the same size.}\label{tab:sim}
\small
\begin{tabular}{lccccc}
\toprule
$N$ & $\xi(\hat C_N)$ & $\xi_N$ & $\Rrho(\hat C_N)$ & $\tau_{0.25,0.25}(\hat C_N)$ & $\tau_{0.1,0.25}(\hat C_N)$\\
\midrule
\multicolumn{6}{l}{\emph{Gaussian}}\\
true value & 0.114 & 0.114 & 0.431 & 0.307 & 0.316\\
512 & 0.115 (0.019) & 0.113 (0.030) & 0.421 (0.035) & 0.269 (0.044) & 0.227 (0.045)\\
1728 & 0.115 (0.010) & 0.113 (0.017) & 0.427 (0.019) & 0.289 (0.027) & 0.272 (0.034)\\
8000 & 0.114 (0.005) & 0.115 (0.008) & 0.429 (0.009) & 0.300 (0.013) & 0.297 (0.017)\\
46656 & 0.114 (0.002) & 0.115 (0.003) & 0.430 (0.004) & 0.304 (0.006) & 0.309 (0.008)\\
125000 & 0.114 (0.001) & 0.114 (0.002) & 0.430 (0.002) & 0.305 (0.004) & 0.312 (0.005)\\
\addlinespace
\multicolumn{6}{l}{\emph{Student $t_4$}}\\
true value & 0.114 & 0.114 & 0.414 & 0.362 & 0.403\\
512 & 0.114 (0.021) & 0.112 (0.032) & 0.409 (0.038) & 0.307 (0.048) & 0.270 (0.049)\\
1728 & 0.114 (0.011) & 0.115 (0.018) & 0.413 (0.021) & 0.334 (0.029) & 0.334 (0.037)\\
8000 & 0.114 (0.006) & 0.114 (0.009) & 0.414 (0.010) & 0.350 (0.014) & 0.372 (0.019)\\
46656 & 0.114 (0.002) & 0.115 (0.003) & 0.414 (0.004) & 0.358 (0.006) & 0.392 (0.008)\\
125000 & 0.114 (0.001) & 0.114 (0.002) & 0.414 (0.002) & 0.360 (0.003) & 0.397 (0.005)\\
\addlinespace
\multicolumn{6}{l}{\emph{Resolution bias $\kappa(C_m)-\kappa(C)$, Student $t_4$}}\\
512 & $-$0.011 & -- & $-$0.018 & $-$0.057 & $-$0.133\\
1728 & $-$0.006 & -- & $-$0.009 & $-$0.029 & $-$0.068\\
8000 & $-$0.003 & -- & $-$0.003 & $-$0.012 & $-$0.030\\
46656 & $-$0.001 & -- & $-$0.001 & $-$0.004 & $-$0.010\\
125000 & $-$0.001 & -- & $-$0.001 & $-$0.002 & $-$0.006\\
\addlinespace
\multicolumn{6}{l}{\emph{Probability that the Student estimate exceeds the Gaussian estimate}}\\
512 & 0.47 & 0.49 & 0.40 & 0.72 & 0.73\\
1728 & 0.47 & 0.53 & 0.32 & 0.87 & 0.89\\
8000 & 0.46 & 0.50 & 0.12 & 1.00 & 1.00\\
46656 & 0.46 & 0.51 & 0.00 & 1.00 & 1.00\\
125000 & 0.47 & 0.50 & 0.00 & 1.00 & 1.00\\
\bottomrule
\end{tabular}
\end{table}

First, the global measures do not distinguish the two copulas at any sample size: $\xi(\hat C_N)$ and $\xi_N$ estimate the same value, and $\Rrho(\hat C_N)$ ranks the Gaussian copula higher, as the population values do. The tail-localised measures rank the Student copula higher with probability $0.87$--$0.89$ at $N=1728$ and with probability at least $0.996$ from $N=8000$ on. Second, the root mean squared errors of all estimators decrease like $N^{-1/2}$; least squares fits of $\log(\text{RMSE})$ against $\log N$ over $N\ge1728$ give slopes between $-0.54$ and $-0.49$ for both models. Third, the bias of the tail-localised estimators is almost exactly the resolution bias; it is negative, decreases approximately like $m^{-1.8}$, and dominates the error of $\tau_{0.1,0.25}(\hat C_N)$ for $N\le2000$; the standard deviations of the tail-localised estimators are about one and a half to two times those of $\Rrho(\hat C_N)$. For $\xi(\hat C_N)$ the negative resolution bias is offset by a positive bias due to sampling variability, which is plausible because $\xi$ is a quadratic functional of $\partial_1C$. Fourth, the standard deviation of $\xi(\hat C_N)$ is only about $60\%$ of that of $\xi_N$, in line with the known inefficiency of $\xi_N$ \citep{ShiDrtonHan2022,LinHan2023}; this is a finite-sample observation for two models, and unlike $\xi_N$, for which a central limit theorem is available \citep{LinHan2022}, the asymptotic distribution of the checkerboard estimators is unknown.

In summary, tail-localised measures with moderate tail levels are estimable with useful precision from a few thousand observations and separate copulas that every global measure considers equivalent, whereas extreme tail levels are beyond the resolution of nonparametric estimators. A bias correction based on the deterministic part of the error, limit distributions of the checkerboard estimators, and a data-driven choice of $m$ in the spirit of the cross-validation of \citet{SDS2024} are natural next steps.

\section{Application: systematic risk in large credit portfolios}\label{sec:credit}

We interpret the results of Sections~\ref{sec:dictionary}--\ref{sec:sharp} in the simplest factor model of portfolio credit risk, in which the variables $L_t$ are the objects of primary interest. All facts about the model itself are standard, see \citet{FreyMcNeil2003} and \citet[Chapter~11]{McNeilFreyEmbrechts2015}.

\subsection{The model}

Consider obligors $i=1,2,\dots$ with latent variables $V_i$ that are conditionally independent given a systematic factor $U$, with $(U,V_i)\sim C$ for all $i$. Obligor $i$ defaults within the given horizon if $V_i\le t$, so that $t$ is the common default probability. Conditionally on $U$, defaults are independent with probability $L_t=\partial_1C(U,t)$; this is a one-factor Bernoulli mixture model with mixing variable $L_t$. The copula description of default dependence goes back to \citet{Li2000}. For the Gaussian copula with correlation $r$ one obtains the model of \citet{Vasicek2002} with asset correlation $r^2$. With equal exposures and loss given default equal to one, the loss fraction of the first $n$ obligors is $L_t^{(n)}=\frac1n\sum_{i\le n}\mathbf 1\{V_i\le t\}$, and by the conditional law of large numbers $L_t^{(n)}\to L_t$ almost surely. The limit $L_t$ is the loss of a large homogeneous portfolio.

For a loss $L$ with values in $[0,1]$ and a level $x\in(0,1)$ we write $\operatorname{ES}_x(L)=\frac1x\int_0^xL^*(s)\,\dd s$ for the expected loss in the worst fraction $x$ of scenarios (expected shortfall at confidence level $1-x$) and $\operatorname{VaR}_x(L)=L^*(x)$ for the corresponding quantile. In the asymptotic single risk factor framework of \citet{Gordy2003}, on which the internal ratings-based approach of the Basel accords is built, capital for unexpected losses is proportional to $\operatorname{VaR}_x(L_t)-t$ with $x=0.1\%$ \citep{BCBS2005}. A tranche with attachment point $q_1$ and detachment point $q_2$ absorbs the loss $(L_t-q_1)_+-(L_t-q_2)_+$.

In this model the two diagrams of Definition~\ref{def:diagrams} are the two standard descriptions of portfolio risk: the Lorenz diagram $x\mapsto S_t(x)=x\operatorname{ES}_x(L_t)$ records expected shortfall at all levels, and the stop-loss diagram $q\mapsto\Theta_t(q)$ records the expected losses of all senior tranches. The rearranged copula $\Cr$ describes the model in which, for every default probability $t$, the loss has the same law as under $C$ but is a nonincreasing function of the factor, so that small values of the factor are bad for all default probabilities simultaneously; the Schur order compares models by the convex order of their portfolio losses.

\subsection{Translation of the results}

\begin{corollary}\label{cor:credit}
In the model above, the following holds for every copula $C$.
\begin{enumerate}
\item[\textup{(a)}] \emph{Diversification.} $\Var(L_t^{(n)})=\Var(L_t)+\frac1n\bigl(t(1-t)-\Var(L_t)\bigr)$. The indicator correlation $\varrho_C(t)$ is the default correlation of two obligors and the fraction of the default variance $t(1-t)$ that cannot be diversified away, and $\xi(C)=\int_0^1\varrho_C(t)\,6t(1-t)\,\dd t$.
\item[\textup{(b)}] \emph{Risk measures.} $\operatorname{ES}_x(L_t)=\Cr(x,t)/x$, and $\operatorname{VaR}_x(L_t)=\partial_1\Cr(x,t)$ for almost every $x$. The expected loss of the tranche $[q_1,q_2]$ is $\Theta_t(q_1)-\Theta_t(q_2)$. Moreover $C\le_SD$ holds if and only if, for every $t$, every law-invariant convex monetary risk measure on $L^\infty$, applied to the loss $L_t$ with the convention that larger losses are riskier, is at least as large under $D$ as under $C$.
\item[\textup{(c)}] \emph{Spectral measures.} $N_\phi\,\kappa_\phi(C)$ is the average over all default probabilities $t$ of the spectral risk capital $\mathcal M^+_\phi(L_t)-t$ of the portfolio, and $\kappa_\phi(C)$ is its ratio to the largest possible value. In particular $\tfrac16\Rrho(C)$ is the average of $\E\max\{L_t,L_t'\}-t$, and $\tau_{p,b}(C)$ is the normalised average, over default probabilities $t\le b$, of a spectral risk measure that only takes the worst fraction $p$ of scenarios into account.
\item[\textup{(d)}] \emph{Limits of a single parameter.} Let $x\le t$. There is a copula with $\operatorname{ES}_x(L_t)=1$ and $\xi(C)=2x(1-t)^3/(1-x)^2$, and none with a smaller value of $\xi$; the same holds for $\Rrho$ with the value $3x(1-t)^2/(1-x)$. Every senior tranche $[q,1]$ of a portfolio with default probability $t$ carries its maximal expected loss $(1-q)\,t$ under the block copula $B_t$, for which $\xi=2t(1-t)$, and under no copula with a smaller value of $\xi$.
\item[\textup{(e)}] \emph{Tail-localised control.} For $p\le b$ and $t\le b$, the excess $e=\operatorname{ES}_p(L_t)-t$ satisfies \eqref{eq:stress}; in particular $\operatorname{VaR}_p(L_t)-t\le\operatorname{ES}_p(L_t)-t\le\sqrt{\tfrac13(1-p)(3-2p)\,\tau_p(C)}$.
\item[\textup{(f)}] \emph{Control by default correlation.} $\operatorname{ES}_x(L_t)-t\le\sqrt{\varrho_C(t)\,t(1-t)(1-x)/x}$, and the systematic part $\theta_t(q)$ of the expected loss of the tranche $[q,1]$ satisfies \eqref{eq:profTranche}.
\end{enumerate}
\end{corollary}

\begin{proof}
(a) Given $U$, the variable $nL_t^{(n)}$ is binomial with parameters $n$ and $L_t$, so $\Var(L_t^{(n)})=\Var(L_t)+\frac1n\E[L_t(1-L_t)]$; the rest is Theorem~\ref{thm:dictionary}(i). (b) The first statements are Lemma~\ref{lem:basic}(a) and Definition~\ref{def:diagrams}. Law-invariant convex risk measures on $L^\infty$ of an atomless probability space have the Fatou property \citep{JouiniSchachermayerTouzi2006} and are therefore monotone in the convex order \citep[Chapter~4]{FoellmerSchied2016}, and expected shortfall at all levels characterises the convex order for equal means, so the last statement follows from \eqref{eq:schur}. (c) is Definition~\ref{def:spectral} together with Theorem~\ref{thm:mirror}(a) and Proposition~\ref{prop:two}; (d) is Corollary~\ref{cor:tail} and Theorem~\ref{thm:stoploss}(b); (e) is Corollary~\ref{cor:stress} together with $\operatorname{VaR}_p\le\operatorname{ES}_p$; (f) is Proposition~\ref{prop:profile}.
\end{proof}

Part (b) says that the Schur order is the ordering of systematic risk on which all users of law-invariant convex risk measures agree; this is the order-theoretic content of the comparison results of \citet{AnsariLuetkebohmert2025}. The measure $\xi$ is the aggregate non-diversifiable share of default variance, but it weights all default probabilities and all scenarios evenly; a value such as $\xi=0.002$ is compatible with the default of the entire portfolio in the worst $0.1\%$ of scenarios at a default probability of $1\%$. The extremal copula $T_{x,t,x(1-t)}$ has a simple interpretation: with probability $x$ a crisis state occurs in which every obligor with default probability $t$ defaults, and outside this state the factor carries little information. That models with an identical value of a dependence parameter can differ widely in the tail is well known in this area \citep{FreyMcNeil2003,DonnellyEmbrechts2010}; parts (d)--(f) quantify this for the measures considered here, and part (e) shows that a spectral measure adapted to the tail level and to the relevant range of default probabilities does bound tail risk.

\subsection{Numerical example}

We return to the two copulas of Section~\ref{sec:illustration}. The Gaussian copula with $r^2=0.2$ is Vasicek's model with an asset correlation of $20\%$, a value within the range used in the Basel risk weight functions for corporate exposures \citep{BCBS2005}. The Student copula has the same value $\xi=0.1145$. It should be read as a factor copula linking the factor and the latent variable; it is not the multivariate $t$ threshold model, in which obligors are not conditionally independent given a single factor.

\begin{table}[htbp]
\centering
\caption{Quantile of the loss at level $x=0.1\%$ and default correlation, for three default probabilities $t$. The corresponding values of $\operatorname{ES}_x(L_t)=\Cr(x,t)/x$ are given in Table~\ref{tab:local}.}\label{tab:var}
\begin{tabular}{lcccccc}
\toprule
 & \multicolumn{3}{c}{$\operatorname{VaR}_x(L_t)$} & \multicolumn{3}{c}{$\varrho_C(t)$}\\
\cmidrule(lr){2-4}\cmidrule(lr){5-7}
$t$ & $0.1\%$ & $1\%$ & $5\%$ & $0.1\%$ & $1\%$ & $5\%$\\
\midrule
Gaussian & 0.028 & 0.146 & 0.384 & 0.006 & 0.024 & 0.058\\
Student $t_4$ & 0.119 & 0.427 & 0.628 & 0.082 & 0.093 & 0.105\\
\bottomrule
\end{tabular}
\end{table}

\begin{table}[htbp]
\centering
\caption{Expected loss $\Theta_t(q)$ of the senior tranche $[q,1]$ at default probability $t=1\%$, in basis points of the portfolio notional. The last two rows are the bound \eqref{eq:profTranche} evaluated with the default correlation of the Gaussian model, and the maximal value $t(1-q)$, which is attained by the block copula $B_t$ with $\xi=0.0198$.}\label{tab:tranche}
\begin{tabular}{lcccc}
\toprule
Attachment point $q$ & $3\%$ & $7\%$ & $15\%$ & $30\%$\\
\midrule
Gaussian & 16.7 & 3.6 & 0.32 & 0.007\\
Student $t_4$ & 35.4 & 21.9 & 10.7 & 3.1\\
Bound \eqref{eq:profTranche}, Gaussian $\varrho_C$ & 26.4 & 9.8 & 4.3 & 2.1\\
Maximum $t(1-q)$ & 97 & 93 & 85 & 70\\
\bottomrule
\end{tabular}
\end{table}

\begin{figure}[htbp]
\centering
\includegraphics[width=0.95\textwidth]{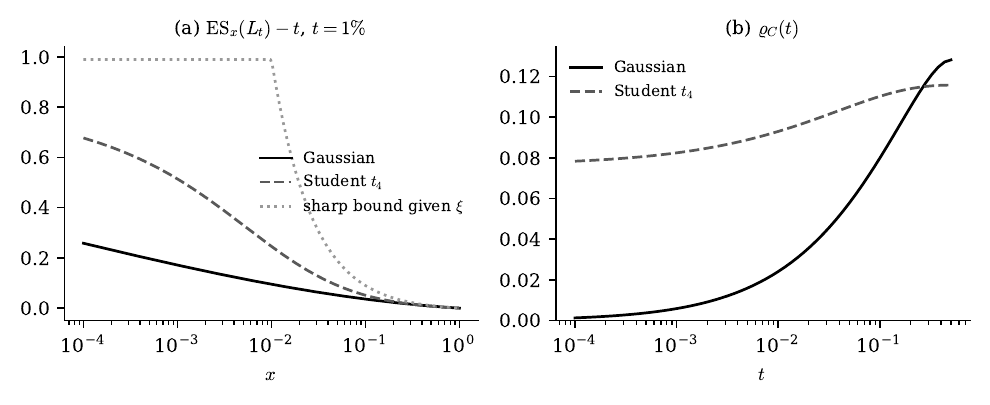}
\caption{(a) The excess $x\mapsto\operatorname{ES}_x(L_t)-t$ at default probability $t=1\%$ for the Gaussian and the Student model with identical $\xi=0.1145$, together with the sharp upper bound given this value of $\xi$, which coincides with the Fr\'echet--Hoeffding bound $M(x,t)/x-t$ at every level $x$. (b) Default correlation profiles $t\mapsto\varrho_C(t)$.}\label{fig:stress}
\end{figure}

Although $\xi$, $\Rrho$ and $\zeta_1$ practically coincide for the two models (Table~\ref{tab:global}), their tail risk differs substantially. At a default probability of $1\%$ the expected loss in the worst $0.1\%$ of scenarios is $18\%$ of the portfolio in the Gaussian and $52\%$ in the Student model (Table~\ref{tab:local}), and the quantile that enters the capital formula is $14.6\%$ against $42.7\%$ (Table~\ref{tab:var}). The default correlation at $t=0.1\%$ differs by a factor of fourteen, in line with the absence of tail dependence in the Gaussian case (Figure~\ref{fig:stress}(b)). For tranches the differences are larger still (Table~\ref{tab:tranche}): the expected loss of a senior tranche attaching at $15\%$ is about thirty times, and that of a tranche attaching at $30\%$ about five hundred times, larger in the Student model. By Figure~\ref{fig:stress}(a), the sharp bound on expected shortfall implied by the common value of $\xi$ is trivial at every level $x$: the Fr\'echet--Hoeffding bound is attained at $(x,t)$ by a tent copula $T$ whose value of $\xi$ is at most $2t(1-t)=0.0198$, and the copulas $(1-\sigma)T+\sigma M$, $\sigma\in[0,1]$, have the same value at $(x,t)$ and take every value of $\xi$ between $\xi(T)$ and $1$. The tail measure $\tau_{0.1\%,5\%}$, in contrast, gives the nontrivial bounds reported in Table~\ref{tab:local}; the bound \eqref{eq:profES} from the default correlation at $t=1\%$ alone gives $0.49$ and $0.96$.

\section*{Statement on AI use}

The author used the AI assistant Claude (Anthropic) in the preparation of this paper: for literature searches, proofreading, proof verification, figure generation, and for writing the code for the numerical examples and the simulation study. All results and proofs have been checked by the author, who takes full responsibility for the content of the paper.

\bibliographystyle{apalike}
\bibliography{references}

\end{document}